\documentclass[a4, 10pt]{amsart}
\usepackage{amssymb,amsmath}
\newtheorem{thm}{\sc Theorem}[section]
\newtheorem{prp}{\sc Proposition}[section]
\newtheorem{lem}{\sc Lemma}[section]
\newtheorem{rem}{\sc Remark}[section]
\newtheorem{cor}{\sc Corollary}[section]
\newtheorem{ex}{\sc Example}[section]
\def\beweis{{\it Proof.~}}

\def\Res{\mathop{\rm res}}
\def\N{\mathbb{N}}
\def\R{\mathbb{R}}

\def\C{\mathbb{C}}
\def\Z{\mathbb{Z}}

\def\P{\mathcal P}

\def\CL{\mathfrak C}

\def\Re{{\rm Re} \,}
\def\Im{{\rm Im} \,}

\def\mod{{\rm mod\,}}

\def\dist{{\rm dist}}

\def\ds{\displaystyle}
\def\ts{\textstyle}

\def\ende{\blacksquare}
 \def\Ende{$\ende$\medskip}

\def\x{\mathfrak{x}}
\def\y{\mathfrak{y}}

\def\z{\mathfrak{z}}
\def\w{\mathfrak{w}}
\def\a{\mathfrak{a}}

\def\c{\mathfrak{c}}

\def\f{\mathfrak{f}}

\def\h{\mathbf{h}}

\def\p{\mathfrak{p}}
\def\q{\mathfrak{q}}

\def\0{\mathfrak{0}}
\def\1{\mathbf{1}}

\def\Y{\mathfrak{Y}}
\def\barY{\overline\Y}\def\barY{\tilde\Y}
\def\doppel{\star\ast}

\def\be#1{\begin{equation}\label{#1}}
\def\ee{\end{equation}}
\def\ko#1{o(|z|^{#1})}
\def\GO#1{O(|z|^{#1})}
\def\schema#1#2#3#4{\left[#3\!\!\begin{array}{c}#2\cr#4\end{array}\!\!#1\right]}

\begin{document}

\begin{center}{\huge A unified approach to the Painlev\'e Transcendents}\bigskip

by\medskip

{\sc Norbert Steinmetz}\medskip

{\small\it Technische Universit\"at Dortmund}
\end{center}

\bigskip\bigskip

{\small\begin{center} {\bf Abstract}\end{center}
\begin{quote} We utilise a recent approach via the so-called re-scaling method to derive a unified and comprehensive theory of the solutions to
{\sc Painlev\'e}'s differential equations (I), (II) and (IV), with emphasis on the most elaborate equation (IV).
\end{quote}


}\bigskip\bigskip

{\small\vspace{5mm}\begin{center}\begin{tabular}{rlr}
\ref{INTRO}.& Introduction&\pageref{INTRO}\\
\ref{VALUED}.& Value Distribution&\pageref{VALUED}\\
\ref{PTYF}.& {Painlev\'e} Transcendents and {Yosida} Functions&\pageref{PTYF}\\
\ref{PFSaAE}.& `Pole-free' Sectors&\pageref{PFSaAE}\\
\ref{ASYMPEXPANS}.& Asymptotic Expansions&\pageref{ASYMPEXPANS}\\
\ref{MULTTRUNC}.& Truncated {Painlev\'e} Transcendents&\pageref{MULTTRUNC}\\
\ref{DICHOTOMY}.& The Dichotomy of the Order&\pageref{DICHOTOMY}\\
\ref{SUBNORMALSOLUTIONS}.& Sub-normal Solutions&\pageref{SUBNORMALSOLUTIONS}\\
\ref{PolesZeros}.& The Distribution of Zeros and Poles&\pageref{PolesZeros}\\
\ref{VALUEDISTRIBUTION}.& Deficient Values and Functions&\pageref{VALUEDISTRIBUTION}\\
\ref{AppendixB}.& Appendix A: The Phragm\'en-Lindel\"of Principle&\pageref{AppendixB}\\
\ref{AppendixC}.&Appendix B: Asymptotic Expansions&\pageref{AppendixC}\\
\end{tabular}\end{center}}\bigskip\bigskip

\section{\bf Introduction}\label{INTRO}
Generally speaking, there are two different approaches to {\sc Painlev\'e}'s differential equations.
The first one is based on the {\sc Riemann-Hilbert} method and the method of isomonodromic deformations, and is strongly linked to the fields of {\it Special Functions} and {\it Mathematical Physics}, while the second one may be viewed as a part of {\it Complex Analysis}.
We exemplarily mention the monographs by {\sc Fokas}, {\sc Its}, {\sc Kapaev}, and {\sc Novoksh\"{e}nov} \cite{FIKN} on one hand, and {\sc Gromak}, {\sc Laine}, and {\sc Shimomura}~\cite{GLS} on the other. Although there is some overlap, the methods, issues, and results are quite different. It is even difficult to translate the results from one language into the other.
The aim of this paper is to develop a unified theory for the solutions to the {\sc Painlev\'e} equations
$$\begin{array}{rrcl}{\rm (I)}\qquad&w''&=&z+6w^2\cr
{\rm (II)}\qquad&w''&=&\alpha+zw+2w^3\cr
{\rm (IV)}\qquad&2ww''&=&w'^2+3w^4+8zw^3+4(z^2-\alpha)w^2+2\beta\end{array}$$
with complex analytic methods. We shall investigate all equations simultaneously and in the same way,  particular attention, however, will be payed
to equation (IV), which usually is neglected in the literature.
Several of the results in this paper have been proved in the past decade by different authors with different methods, in a bid to verify statements of the `mathematical folklore',
many of them dating back to {\sc Boutroux}'s papers \cite{Boutroux1}(\footnote{There is some dissent about the validity of {\sc Boutroux}'s methods.}).
Not every result stated and proved in the present paper is new. What is new is the methodical unification, which is based on the re-scaling method
developed in \cite{NStPainleveOrder} for equation (I), which itself was inspired by the so-called
{\sc Zalcman} {\it re-scaling method} \cite{LZ1, LZ2}.

It is taken for granted that every solution is meromorphic on the plane. For recent proofs the
reader is referred to {\sc Hinkkanen} and {\sc Laine}~\cite{HL1}, {\sc Shimomura}~\cite{ShimomuraA},
and the author~\cite{NStPainleveEx}.
To the convenience of the reader we  defer two major tools
to Appendix A ({\sc Phragm\'en}-{\sc Lindel\"of} Principle) and Appendix B (existence of asymptotic expansions) at the end of the paper.

Any transcendental solution to some {\sc Painlev\'e} equation is called
{\it {\sc Painlev\'e} transcendent}(\footnote{The original meaning of {\it transcendent} was different,
namely: {\it The solutions are transcendental functions of the `two constants of integration'.}}), more precisely,  {\it first, second,} and {\it fourth} transcendent,
respectively. Every {\sc Painlev\'e} equation has a {\it first integral} $W$:
\begin{equation}\label{FirstIntegral}
\left.\begin{array}{lll}
{\rm (i)}& w'^2=4w^3+2zw -2W,&W'=w\cr
{\rm (ii)} &w'^2= w^4+zw^2+2\alpha w-W,&W'=w^2\cr
{\rm (iv)} &w'^2= w^4+4zw^3+4(z^2-\alpha)w^2-2\beta-4wW,\quad& W'=w^2+2zw\cr
\end{array}\right.
 \end{equation}
At every pole $p$, the {\sc Painlev\'e} transcendents and their first integrals have {\sc Laurent} series developments
\begin{equation}\label{wLaurent}
\left.\begin{array}{lrl}
{\rm (i)}&\!\!\!w=&\ds ~(z-p)^{-2}-\frac p{10}(z-p)^2-\frac 16(z-p)^3+\h(z-p)^4+\cdots \cr
 &\!\!\!W=&\ds\!\!\! -(z-p)^{-1}-14\h-\frac p{30}(z-p)^3-\frac 1{24}(z-p)^4+\cdots \cr
{\rm (ii)}&\!\!\! w=&\ds\!\!\!\epsilon (z-p)^{-1}-\epsilon\frac p{6}(z-p)-\frac {\alpha+\epsilon}4(z-p)^2
+\h(z-p)^3+\cdots\cr
&\!\!\!W=&\ds\!\!\!-(z-p)^{-1}+10\epsilon \h-\frac7{36}p^2-\frac p{3}(z-p)
-\frac {1+\epsilon\alpha}4(z-p)^2+\cdots\cr
{\rm (iv)}&\!\!\! w=&\ds\!\!\!\epsilon (z-p)^{-1}-p
+\frac\epsilon 3(p^2+2\alpha-4\epsilon)(z-p)+\h(z-p)^2+\cdots\cr
&\!\!\!W=&\ds\!\!\!-(z-p)^{-1}+2\h+2(\alpha -\epsilon)p+\frac13(4\alpha-p^2-2\epsilon)(z-p)+\cdots
\end{array}\right.
 \end{equation}
respectively ($\epsilon=\pm 1$); the coefficient $\h=\h(p)$ remains undetermined,
and {\it free}: the pole $p$, the sign $\epsilon$, and $\h$
may be prescribed to define a unique solution in the same way than
initial values $w_0$ and $w'_0$ at $z_0$ do. The significance of $\h$ cannot be overestimated.

\section{\bf Value Distribution}\label{VALUED}

According to {\sc Boutroux}~\cite{Boutroux1} the {\sc Painlev\'e} transcendents
have order of growth at most $5/2$, $3$, and $4$, respectively. These estimates were confirmed
independently and with different methods by {\sc Shimomura}~\cite{Shimomura1} and the author~\cite{NStPainleveOrder}.
In any case we have $m(r,w)=O(\log r)$, thus the growth of the {\sc Nevanlinna} characteristic $T(r,w)$ is governed by the counting functions of poles $N(r,w)$ and $n(r,w)$.
For notation and results in {\sc Nevanlinna} theory the reader is referred to {\sc Hayman}'s monograph~\cite{Hayman}.

\subsection{\it The key estimates}Analysing section 5.2 in {\sc Shimomura}'s paper \cite{Shimomura1} or section~6 in the author's paper \cite{NStPainleveOrder} yields the following

\begin{prp} Let $w$ be any fourth transcendent, denote by $\P$ the set of non-zero poles of $w$, and set
\be{TRI}\triangle_\delta(p)=\{z:|z-p|<\delta|p|^{-1}\}\quad{\rm and}\quad \P_\delta={\ts\bigcup_{p\in\P}}\triangle_\delta(p).\ee
Then for $\delta>0$ sufficiently small, the discs $\triangle_\delta(p)$ are mutually disjoint and
\be{KEY}w=\GO{}\quad(z\to\infty, ~z\notin\P_\delta).\ee\end{prp}

\medskip This leads easily to the following estimates.

\begin{prp}\label{PROP1IV}For any fourth transcendent it is true that
\be{wWestimate}w'=\GO{2} {\rm ~and}~W=\GO{3} \quad(z\to\infty,~z\notin\P_\delta),\ee
\be{YBed}w'=O(|z|^2+|w|^2)\quad(z\to\infty~{\rm without~restriction}).\ee\end{prp}

\beweis The estimate of $w'$ and $W$ in (\ref{wWestimate}) follows from {\sc Cauchy}'s integral theorem for the first derivative applied to the circle $|\zeta-z|=\frac\delta 2|z|^{-1}$, and $W'=w^2+2zw$, respectively, on combination with (\ref{KEY}). To prove the last assertion, we note that $f(z)=\epsilon w'(z)+w(z)^2+zw(z)$ is regular at any pole $p$ with residue $\epsilon$, and satisfies $f(z)=\GO{2}$
on $\partial\triangle_\delta(p)$, hence also on $\triangle_\delta(p)$ by the maximum principle. This yields
$|w'|\le\GO{2}+|w|^2+|zw|=\GO{2}+O(|w|^2).$ \Ende

\begin{rem}\rm The data for the first and second transcendents reads as follows:
$$\triangle_\delta(p)=\{z:|z-p|<\delta|p|^{-\frac 14}\},\leqno{(\ref{TRI}')}$$
$$w=\GO{\frac 12}\quad(z\to\infty,~z\notin\P_\delta),\leqno{(\ref{KEY}')}$$
$$w'=\GO{\frac 34} {\rm ~and}~W=\GO{\frac 32} \quad(z\to\infty,~z\notin\P_\delta),\leqno{(\ref{wWestimate}')}$$
$$w'=O(|z|^{\frac 34}+|w|^{\frac 32})\quad(z\to\infty~{\rm without~restriction}),\leqno(\ref{YBed}')$$
and
$$\triangle_\delta(p)=\{z:|z-p|<\delta|p|^{-\frac 12}\},\leqno{(\ref{TRI}'')}$$
$$w=\GO{\frac 12}\quad(z\to\infty,~z\notin\P_\delta),\leqno{(\ref{KEY}'')}$$
$$w'=\GO{},{\rm ~and} W=\GO{2} \quad(z\to\infty,~z\notin\P_\delta),\leqno{(\ref{wWestimate}'')}$$
$$w'=O(|z|+|w|^{2})\quad(z\to\infty~{\rm without~restriction}).\leqno(\ref{YBed}'')$$\end{rem}

\subsection{\it The spherical derivative}The value distribution of any meromorphic function takes place in regions
where the spherical derivative
$$f^\sharp(z)=\frac{|f'(z)|}{1+|f(z)|^2}$$
is large, while $f$ behaves tame where $f^\sharp$ is small.

\begin{prp}\label{PROP4IV}Let $w$ be any fourth transcendent. Then
$f(z)=w(z)/z$ has spherical derivative
$$f^\sharp(z)=\GO{}.$$
Similarly, for any first and second transcendent, $f(z)=w(z)^2/z$ has spherical derivative
$f^\sharp(z)=\GO{\frac 14}$ and $f^\sharp(z)=\GO{\frac 12},$
respectively.\end{prp}

\beweis In the first case we obtain from (\ref{YBed})
$$f^\sharp(z)\le\frac{|zw'|}{|z|^2+|w|^2}+\frac{|w|}{|z|^2+|w|^2}=\GO{}+O(1).$$
In the same manner ($\ref{YBed}'$) gives
$$\begin{array}{rcl}
\ds\frac{|f'(z)|}{1+|f(z)|^{\frac 32}}&\le&\ds\frac{2|ww'|}{|z|(1+|w^2/z|^{\frac 32})}+\frac{|w|^2}{|z|^2(1+|w^2/z|^{\frac 32})}\cr
&=&\ds|z|^{\frac14}\frac{|w'|}{|z|^{\frac 34}+|w|^{\frac 32}}\frac{2|z|^{\frac14}|w|(|z|^{\frac 34}+|w|^{\frac 32})}{|z|^{\frac 32}+|w|^3}+O(1)=\GO{\frac 14}\end{array}$$
(the term $\ds\frac{2xy^2(x^{3}+y^{3})}{x^6+y^6}=\frac{2t^2(1+t^3)}{1+t^6}$ with $x=|z|^{\frac14}>0$, $y=|w|^{\frac 12}>0$ and $t=y/x>0$ is bounded), {\it a fortiori} $f^\sharp(z)=\GO{\frac 14}$.
Similarly, ($\ref{YBed}''$) yields
$$\begin{array}{rcl}
f^\sharp(z)&\le&\ds\frac{2|ww'|}{|z|(1+|w^2/z|^2)}+\frac{|w|^2}{|z|^2(1+|w^2/z|^2)}\cr
&=&\ds|z|^{\frac 12}\frac{|w'|}{|z|+|w|^2}\frac{2|z|^{\frac 12}|w|(|z|+|w|^2)}{|z|^2+|w|^4}+O(1)=\GO{\frac 12}.~\ende\cr
\end{array}$$

\subsection{\it The order of growth}\label{GAMMA}There are several possibilities to prove
\be{ORDER}T(r,w)=O(r^{\frac 52}),~T(r,w)=O(r^3),~{\rm and~} T(r,w)=O(r^4)\ee
for first, second, and fourth transcendents, respectively:
\begin{itemize}\item[{\bf a.}] {\it Polar statistics.} The fact that the discs $\triangle_\delta(p)$ are mutually disjoint and have area $\pi\delta^2|p|^{-1/2}$, $\pi\delta^2|p|^{-1}$, and $\pi\delta^2|p|^{-2}$ implies $n(r,w)r^{-\frac12}=O(r^2)$, $n(r,w)r^{-1}=O(r^2)$, and
$n(r,w)r^{-2}=O(r^2)$, respectively, hence the assertion follows from $m(r,w)=O(\log r)$ in any case.
\item[{\bf b.}] {\it Residue Theorem.} It is not hard to construct a curve $\Gamma_r$ of length $O(r)$ which encloses exactly the poles of $w$ with $|p|\le r$
and such that $W=\GO{\frac 32}$,  $W=\GO{2}$, and $W=\GO{3}$, respectively, holds on $\Gamma_r$: starting with the circle $C_r:|z|=r$, we replace any sub-arc $C_r\cap\triangle_\delta(p)$ by a sub-arc of $\partial \triangle_\delta(p)$---the part outside $|z|=r$ if $|p|\le r$, and the part inside $|z|=r$ otherwise.
Then the estimate for $n(r,w)$ given in {\bf a.}\ follows from
$$n(r,w)=-\frac1{2\pi i}\int_{\Gamma_r}W(z)\, dz\quad(\Res_p W=-1).$$
\item[{\bf c.}] {\it The {\sc Ahlfors}-{\sc Shimizu} characteristic}  of $f$ is given by
$$T(r,f)=\frac1\pi\int_0^r A(t,f)\,\frac{dt}t\quad{\rm with}~ A(t,f)=\int_{|z|<t}f^\sharp(z)^2\,d(x,y).$$
From Proposition~\ref{PROP4IV} it follows that $A(t,f)=O(t^{2\lambda+2})$ holds, with $\lambda=\frac14,$ $\lambda=\frac 12,$ and $\lambda=1$, respectively, hence
also $T(r,f)=O(r^{2\lambda+2})$ and $T(r,w)=O(r^{2\lambda+2})$. We note, however, that it is impossible to derive sharp bounds for $T(r,w)$ solely from the
sharp estimates $w^\sharp=\GO{\frac 34}$, $w^\sharp=\GO{\frac 32}$, and $w^\sharp=\GO{2}$, respectively, which just lead to $T(r,w)=O(r^{\frac 72})$, $T(r,w)=O(r^5)$, and $T(r,w)=O(r^6)$. The reason for this is that $w^\sharp(z)$ is much too large on small  neighbourhoods of the zeros of $w$, while the density ${f^\sharp}(z)^2$ of $A(t,f)$ is uniformly distributed but not too large.
\item[{\bf d.}] For any transcendent, the {\it entire} function $\ds F(z)=e^{-\int W(z)\,dz}$
has simple zeros at the poles of $w$  and satisfies
$\log^+|F(z)|=\GO{\varrho}$ as $z\to\infty$ outside $\P_\delta$ $(\varrho=\frac 52,3,4$, respectively),
and this also holds inside the discs
$\triangle_\delta(p)$ by the maximum principle. This implies $T(r,F)=m(r,F)=O(r^\varrho)$ and
$T(r,W)=T(r,F'/F)=O(r^\varrho).$
\end{itemize}

\begin{rem}\rm  In accordance with {\sc Shimomura} (\cite{Shimomura1}, p.\ 259) we note that the proofs of the estimates for $n(r,w)$ in {\sc Hille}~\cite{EH}, p.~443,
and also in {\sc Kitaev}~\cite{KITAEV}, p.~134 (with reference to {\sc Hille}), are incorrect;
both `proofs' implicitly make use of the key fact that  $w(z)^{-3}W(z)$ in case (I) and  $w(z)^{-4}W(z)$ in case (II) (denoted $J(z)$ in \cite{EH, KITAEV}) is bounded
on $\{z:|w(z)|>|z|^{1/2}\}$, which is more or less equivalent to what has to be proved. Several other attempts like \cite{HWS, HW}, working solely with $w^\sharp$, failed
by the reason outlined in {\bf c.}\end{rem}

\section{\bf Painlev\'e Transcendents and Yosida Functions}\label{PTYF}
\subsection{\it Re-scaling}Let $a$ and $b>-1$ be real parameters. The class $\barY_{a,b}$ consists of all meromorphic functions $f$ such that
the family $(f_h)_{|h|>1}$ of functions
\be{Resc}f_h(\z)=h^{-a}f(h+h^{-b}\z)\ee
is normal on $\C$ in the sense of {\sc Montel}, and all limit functions  $\f=\lim_{h_n\to\infty}f_{h_n}$ are $\not\equiv\infty$, at least one of them being non-constant.
If, in addition, {\it all} limit functions are non-constant, then $f$ is said to belong to the {\it {\sc Yosida} class}
$\Y_{a,b}$. The functions of class $\Y_{0,0}$ were introduced by {\sc Yosida}~\cite{Yosida2}, and for arbitrary real parameters by the author~\cite{NStYosida}.
The class $\Y_{0,0}$ is universal in the sense that it contains all limit functions $\f=\lim_{h_n\to\infty}f_{h_n}$
for $f\in\Y_{a,b}$. The functions $f\in\barY_{a,b}$  have
striking properties, for example they satisfy $f^\sharp(z)=\GO{|a|+b}$, $T(r,f)=O(r^{2+2b})$ and $m(r,f)=O(\log r)$, and even $T(r,f)\asymp r^{2+2b}$ and $m(r,1/f')=O(\log r)$ if $f\in\Y_{a,b}$.

\begin{lem}Normality of {\it any} re-scaled family $(f_h)_{|h|>1}$ is equivalent to
\be{NH}\limsup_{z\to\infty}\frac{|f'(z)||z|^{a-b}}{|z|^{2a}+|f(z)|^2}<\infty,\ee
while $f\in\Y_{a,b}$ is equivalent to
\be{NHa}\ds\liminf_{h\to\infty}\sup_{|z-h|<\delta|h|^{-b}}\frac{|f'(z)||z|^{a-b}}{|z|^{2a}+|f(z)|^2}>0\ee
for some (every) $\delta>0$, together with (\ref{NH}).\end{lem}

\beweis To prove the necessity of (\ref{NH}) we just note that $\ds f_h^\sharp(0)=\frac{|h|^{-b-a}|f'(h)|}{1+|h|^{-2a}|f(h)|^2}$ has to be bounded on $|h|>1$ by {\sc Marty}'s criterion.
Conversely, if $R>0$ is arbitrary, $|\z|<R$ and $z=h+h^{-b}\z=h(1+O(|h|^{-1-b}))$, we obtain
$$\limsup_{h\to\infty}f_h^\sharp(\z)=\limsup_{h\to\infty}\frac{|h|^{-b-a}|f'(z)|}{1+|h|^{-2a}|f(z)|^2}=\limsup_{z\to\infty}\frac{|z|^{-b-a}|f'(z)|}{1+|z|^{-2a}|f(z)|^2},$$
uniformly with respect to $\z$, hence $f^\sharp_h(\z)$ is bounded on $|\z|<R$, $|h|>1$, and normality of $(f_h)_{|h|>1}$ follows from  {\sc Marty}'s criterion.
Finally, condition (\ref{NHa}) is equivalent to the fact that every limit function $\f=\lim_{h_n\to\infty}f_{h_n}$ is non-constant.~\Ende

We note that by $(\ref{YBed}'), (\ref{YBed}''),$ and $(\ref{YBed})$ the families $(w_h)_{|h|>1}$ of functions
$$h^{-\frac12}w(h+h^{-\frac14}\z),~h^{-\frac12}w(h+h^{-\frac12}\z), {\rm ~and~} h^{-1}w(h+h^{-1}\z)$$
for first, second and fourth {\sc Painlev\'e} transcendents $w$, respectively, are normal.

\begin{thm}\label{PainleveYosida}The first, second, and fourth {\sc Painlev\'e} transcendents belong to the {\sc Yosida} classes
$\barY_{\frac12,\frac14},$ $\barY_{\frac12,\frac12},$ and $\barY_{1,1},$ while the corresponding functions $f$ in Proposition~\ref{PROP4IV}
belong to the classes $\barY_{0,\frac14},$ $\barY_{0,\frac12},$ and $\barY_{0,1},$
respectively. Every limit function $\w=\lim_{h_n\to\infty}w_{h_n}$ satisfies
\be{PainleveElliptic}\begin{array}{lrcl}
{\rm (i)}& \w''&=&6\w^2+1\cr
&\w'^2&=&4\w^3+2\w-2\c\cr
{\rm (ii)}& \w''&=&2\w^3+\w\cr
&\w'^2&=&\w^4+\w^2-\c\cr
{\rm (iv)}& 2\w\w''&=&\w'^2+3\w^4+8\w^3+4\w^2\cr
&\w''&=&2\w^3+6\w^2+4\w-2\c\cr
&\w'^2&=&\w^4+4\w^3+4\w^2-4\c\w,
\end{array}\ee
respectively, with $\c=\c(\{h_n\},w)$ some constant.  \end{thm}

\beweis Normality of the families $(w_h)$ has already been proved. Every limit function $\w=\lim_{h_n\to\infty}w_{h_n}$ is $\not\equiv\infty$ since it has a pole at $\z=0$ if $\lim\limits_{n\to\infty}|h_n|\dist(h_n,\P)=0$, and
is bounded on $|\z|<\delta$ if $\liminf\limits_{n\to\infty}|h_n|\dist(h_n,\P)>0$; in particular it follows that $w$ belongs to the respective {\sc Yosida} class.
In case of equation (IV), $\w$ satisfies $2\w\w''=\w'^2+3\w^4+8\w^3+4\w^2;$
differentiating this equation yields $2\w\w'''=(12\w^2+24\w+8)\w\w',$
hence either $\w=0$ or else $\w'''=(6\w^2+12\w+4)\w'$ and
$\w''=2\w^3+6\w^2+4\w-2\c$
otherwise; on combination with the first equation (\ref{PainleveElliptic})(iv) this leads to
$\w'^2=\w^4+4\w^3+4\w^2-4\c\w.$ For first and second transcendents it is even easier to derive the corresponding differential equations.
Also the assertion on $f$ follows immediately from Proposition~\ref{PROP4IV} and (\ref{NH}). \Ende

\subsection{\it Elliptic and trigonometric limit functions}
Any limit function $\w=\lim_{h_n\to\infty}w_{h_n}$ satisfies some differential equation
\be{ELLI}\w'^2=P(\w;\c),\ee
where $P$ is a polynomial of degree three or four that depends on a complex para\-meter $\c$. Constant limit functions  (not: `constant solutions')
may only occur under special circumstances, namely
if (i)~$6\w^2+1=2\w^3+\w-\c= 0$, (ii)~$2\w^3+\w=\w^4+\w^2-\c=0$, and
(iv)~$3\w^4+8\w^3+4\w^2=\w^3+3\w^2+2\w-\c=0$, hence
\be{Constantsolutions}\begin{array}{lrl}
{\rm (i)}&\w\equiv\pm i/\sqrt{6}&{\rm if~} \c=\pm i\sqrt{2/27}\cr
{\rm (ii)}&\w\equiv 0&{\rm if~} \c=0\cr
&\w\equiv \pm i/\sqrt 2&{\rm if~} \c=-1/4\cr
{\rm (iv)}&\w\equiv 0, -2 &{\rm if~} \c=0\cr
&\w\equiv -2/3 &{\rm if~} \c=-8/27.\end{array}\ee
 In the respective cases (\ref{PainleveElliptic}),  $P$ has discriminant
(i) $2+27\c^2,$  (ii) $\c(4\c+1)^2,$ and (iv) $\c^3(27\c+8),$
and the zeros of the discriminant correspond to the {\it exceptional} parameters noted in (\ref{Constantsolutions}).
It is not hard to integrate equation (\ref{PainleveElliptic}) in the exceptional cases; only (ii) and (iv), and actually only the periods $\omega$ will be of interest.\medskip
\be{table1}
\begin{array}{rlll}{\rm (ii)}&\c=0& \w=\pm 1/\sinh \z& \omega=\pi i\cr
&\c=-1/4&\w=\pm\tan(\z/\sqrt 2)/\sqrt 2&\omega=\pi\sqrt 2 \cr
&&&\cr
{\rm (iv)}&\c=0&\w=2e^{\pm 2\z}/(1-e^{\pm 2\z})&\omega=\pi i\cr
&\c=-8/27&\w=8/(9\tan^2(\z/\sqrt{3})-3)&\omega=\pi\sqrt 3\end{array}\ee

If $P$ has mutually distinct zeros, then the non-constant solutions to (\ref{ELLI}) are elliptic
functions of elliptic order two. The corresponding lattice $\Lambda=\Lambda_\c$ has a basis $(\omega_1,\omega_2)$ such that $|\omega_1|\le|\omega_2|\le|\omega_1\pm\omega_2|$;
$|\omega_1|$ and $|\omega_2|$ are uniquely determined and depend continuously on the parameter $\c$.
Also $\omega_2$ tends to infinity as
$\c$ tends to any exceptional parameter, while $\omega_1$ has a definite limit (namely $\omega$ in (\ref{table1})).

\subsection{\it The cluster set} The constants of integration in the respective equations
(\ref{PainleveElliptic}) are closely related to the corresponding first integrals. Inserting $\w=\lim_{h_n\to\infty} w_{h_n}$ into equation (\ref{FirstIntegral}) we
obtain
\be{LIMIVW}\begin{array}{lll}
{\rm (i)}&\c=\ds\lim_{h_n\to\infty}h_n^{-3/2}W(h_n)&(\inf|h_n|^{1/4}\dist(h_n,\P)>0)\cr
{\rm (ii)}&\c=\ds\lim_{h_n\to\infty}h_n^{-2}W(h_n)&(\inf|h_n|^{1/2}\dist(h_n,\P)>0)\cr
{\rm (iv)}&\c=\ds\lim_{h_n\to\infty}h_n^{-3}W(h_n)&(\inf|h_n|\dist(h_n,\P)>0)\cr\end{array}\ee
respectively. For any  {\sc Painlev\'e} transcendent $w$ we denote by $\CL(w)$ the set of all respective limits (\ref{LIMIVW});
$\CL(w)$ is called the {\it cluster set}
of $w$. We note that in case (IV), say, the family $(W_h)_{|h|>1}$, $W_h(\z)=h^{-3}W(h+h^{-1}\z)$ is {\it quasi-normal} (for the definition see \cite{CTC}):
if $w_{h_n}$ tends to $\w$, then $W_{h_n}$ tends to some constant $\c\in\CL(w)$, locally uniformly on $\C\setminus\{$poles of $\w\}$.

\begin{prp}\label{CLSET}In all cases the cluster set is closed, bounded, and connected, and contains the respective limits
\be{hcluster}\begin{array}{ll}
{\rm (i)}&\ds\lim_{p_n\to\infty}-14\h(p_n)p_n^{-3/2}\cr
{\rm (ii)}&\ds\lim_{p_n\to\infty}10\epsilon_n\h(p_n)p_n^{-2}-7/36\quad(\epsilon_n=\Res_{p_n}w)\cr
{\rm (iv)}&\ds\lim_{p_n\to\infty}2\h(p_n)p_n^{-3},\end{array}\ee
where $(p_n)$ denotes any suitably chosen sequence of poles of $w$. \end{prp}

\beweis Again we will restrict ourselves to equation (IV). For $\delta>0$ sufficiently small, the closed discs $\bar\triangle_\delta(p)$ about the poles $p\ne 0$ are mutually disjoint,
hence the domain $D_\delta=\C\setminus\bigcup_{p\in\P}\bar\triangle_\delta(p)$ is locally path-connected at infinity: any two points $a,b\in D_\delta$ may be joined by a curve contained in $D_\delta$
of spherical length comparable to the spherical distance of $a$ and $b$. We denote the corresponding cluster set of $z^{-3}W(z)$ as $z\to\infty$ on $D_\delta$ by $\CL_\delta(w)$, and
note that $\CL(w)=\bigcup_{\delta>0}\CL_\delta(w)$ and $\CL_{\delta_1}(w)\subset\CL_{\delta_2}(w)$ if $\delta_2<\delta_1$. The cluster sets $\CL_\delta(w)$ are closed, bounded
in $\C$ since $z^{-3}W(z)$ is uniformly bounded on $D_\delta$, and connected by the special property of $D_\delta$.
It remains to show that $\CL_\eta(w)\subset\CL_\delta(w)$ holds for $\delta>\eta$ sufficiently small (such that the discs $\triangle_\delta(p)$ are mutually disjoint), and that the limits (\ref{hcluster})(iv) belong to $\CL_\delta(w)$. If
$\eta\le |h_n|\dist(h_n,\P)<\delta$ and $\dist(h_n,\P)=|h_n-p_n|$
holds for some $p_n\in\P$, and if $\w=\lim_{h_n\to\infty}w_{h_n}$ exists, we replace $h_n$ by $p_n$
with the following effect: from $p_n=h_n+h_n^{-1}\z_n$ with $(\z_n)$ bounded, hence $\z_n\to\z_0$ as we may assume, it follows that
$p_n^{-1}=(1+o(1))h_n^{-1}$, $w_{h_n}(\z)=(1+o(1))w_{p_n}(\z-\z_0+o(1)),$ and
$$\hat\w(\z)=\lim_{p_n\to\infty}w_{p_n}(\z)=\w(\z+\z_0).$$
Thus $\hat\w$ and $\w$ satisfy the the same differential equation, this showing that $\c=\lim_{h_n\to\infty}h_n^{-3}W(h_n)\in\CL_\eta(w)$ coincides with one of the limits (\ref{hcluster})(iv).
Finally, if we start with some sequence $(w_{p_n})$ with $p_n\in\P$, we may as well consider $(w_{h_n})$
with $|h_n||h_n-p_n|=\delta$ without changing the constant of integration, this showing that the limit $\c=\lim_{p_n\to\infty}2\h(p_n)p_n^{-3}$ actually
belongs to $\CL_\delta(w)$. We also obtain $\CL_\delta(w)=\CL_\eta(w)$ by combining both arguments. \Ende

\section{\bf `Pole-free' Sectors}\label{PFSaAE}

\subsection{\it Re-scaling and `pole-free' sectors}
Let $f$ be meromorphic on some sector $S:|\arg z-\theta_0|<\eta$. Then $S$ is called {\it `pole-free'} for $f$, if for every $\delta>0$, $f$ has
only finitely many poles on $S_\delta:|\arg z-\theta_0|<\eta-\delta.$ If $w$ is any fourth {\sc Painlev\'e} transcendent with `pole-free' sector $S$, then the re-scaling process for
sequences $(h_n)$ in $S_\delta$ leads to constant limit functions $\w=\lim_{h_n\to  \infty}w_{h_n}$. More precisely, these constants are contained in $\{0,-2,-2/3\}$ and also in the cluster set
of $w/z$ restricted to $z\to\infty$ on $S$, which again is (compact) and connected. This proves that $w/z$ tends to one of the constants $0$, $-2$, and $-2/3$.
The same argument works also for first and second transcendents, so that in  any `pole-free' sector of any {\sc Painlev\'e} transcendent the following is true:
\be{ASYMPS}\begin{array}{llcl}
{\rm (i)}\phantom{_a}&w=\sqrt{-z/6}+\ko{1/2}&{\rm and}&W(z)=\sqrt{-2/27}z^{3/2}+\ko{3/2}\cr
{\rm (ii)}_a&w=\ko{1/2}&{\rm and}&W(z)=\ko{2}\cr
{\rm (ii)}_b&w=\sqrt{-z/2}+\ko{1/2}&{\rm and}&W(z)=-\frac14 z^2+\ko{2}\cr
{\rm (iv)}_a&w=-\frac23z+\ko{}&{\rm and}&W(z)=-\frac8{27}z^3+\ko{3}\cr
{\rm (iv)}_b&w=-2z+\ko{}&{\rm and}&W(z)=\ko{3}\cr
{\rm (iv)}_c&w=\ko{}&{\rm and}&W(z)=\ko{3}\end{array}\ee
(for some branch of the square-root in (i) and (ii)$_b$) as $z\to\infty$, uniformly on every closed sub-sector of $S$
in the respective cases.

\subsection{\it Stokes sectors} The so-called {\it {\sc Stokes} sectors} and corresponding {\it {\sc Stokes} rays}
$$\begin{array}{lll}
{\rm (i)}&\Sigma_\nu:|\arg z-2\nu\pi/5|<\pi/5&\sigma_\nu:\arg z=(2\nu+1)\pi/5\cr
{\rm (ii)}_a&\Sigma_\nu:|\arg z-2\nu\pi/3|<\pi/3&
\sigma_\nu:\arg z=(2\nu+1)\pi/3\cr
{\rm (ii)}_b&\Sigma_\nu:|\arg z-(2\nu+1)\pi/3|<\pi/3&
\sigma_\nu:\arg z=(2\nu+2)\pi/3\cr
{\rm (iv)}_a&\Sigma_\nu:|\arg z-(2\nu+1)\pi/4|<\pi/4
&\sigma_\nu:\arg z=(\nu+1)\pi/2\cr
{\rm (iv)}_{b,c}&\Sigma_\nu:|\arg z-\nu\pi/2|<\pi/4
&\sigma_\nu:\arg z=(2\nu+1)\pi/4\cr
\end{array}$$
will play an extraordinary role; note that $\sigma_\nu$ separates the adjacent {\sc Stokes} sectors $\Sigma_\nu$ and $\Sigma_{\nu+1}$.

\begin{prp}\label{RAYTOSECTOR}Let $w$ be any {\sc Painlev\'e} transcendent such that some of the corresponding asymptotics {\rm (\ref{ASYMPS})} holds as $z\to\infty$ on some single
non-{\sc Stokes} ray $\arg z=\bar\theta$. Then  the {\sc Stokes} sector which contains that ray is `pole-free' for $w$.\end{prp}

\beweis Again our focus is on fourth transcendents $w$, and for definiteness we assume that the ray $\arg z=\bar\theta$ in question is contained in
$\Sigma_\nu:|\arg z-\nu\pi/2|<\pi/4$, such that $w(z)=-2z+\ko{}$ holds as $z=re^{i\bar\theta}\to\infty$. Then also $W'=2zw+w^2=\ko{2}$ and
$W(z)=\ko{3}$ hold as $z=re^{i\bar\theta}\to\infty$,
and the re-scaling method along any sequence $h_n=r_ne^{i\bar\theta}$ yields $\c=0$ and limit functions satisfying
$$\w'^2=\w^2+4\w^3+4\w^2=\w^2(\w+2)^2,~\w(0)=-2,$$
hence $\w\equiv -2$. Thus for $r\ge r_0$ there exist pole-free discs $|z-re^{i\bar\theta}|<\rho(r)r^{-1}$ such that $\rho(r)\to\infty$ as $r\to\infty$.
We define a sequence $(r_n)$ recursively by $r_{n+1}=r_n+8r_n^{-1}$ and denote by $\theta_n$ $(\bar\theta<\theta_n\le\bar\theta+2\pi)$ the largest angle such that $w$ has no poles on
$$A_n=\{re^{i\theta}:r_n\le r\le r_{n+1},~\bar\theta\le \theta<\theta_n\};$$
$r_0$ is chosen sufficiently large to ensure that $\theta_0<\bar\theta+2\pi$.
Let $J$ denote the set of integers $n$ such that $w$ has at least one pole $z_n$ on $\partial A_n$ (actually $\arg z_n=\theta_n$). Since $w$ is transcendental, $J$ is an
infinite set, and re-scaling about any sequence $(z_{n_\nu})$ with $n_\nu\in J$ then leads to one of the limit functions $\ds\frac{2e^{\pm 2\z}}{1-e^{\pm 2\z}}$
with period $\pi i$ and poles $k\pi i$. This shows that, for $\nu\ge \nu_0$,
\begin{itemize}\item[a.] $w$ has poles $z_{n_\nu}+k(\pi i+o(1))z_{n_\nu}^{-1}$ $(-2\le k\le 2)$
 on $|z-z_{n_\nu}|<8|z_{n_\nu}|^{-1}$, and no others---note that $2\pi<8<3\pi$, and
\item[b.] the poles $z_{n_\nu\pm 1}$ also exist, that is, $n_\nu\pm 1\in J$.\end{itemize}
Thus $J$ contains all integers $n\ge n_0$, and
the sequence $(z_n)_{n\ge n_0}$ is a sub-sequence of some sequence $(p_k)$ of poles satisfying
$$p_{k+1}=p_k+(\omega+o(1)) p_k^{-1}$$
with $\omega\in\{-\pi i, \pi i\}$ independent of $k$. By Lemma~\ref{lemma5} below this implies
$$\theta_n=\arg z_n=\arg p_{k_n}\to (2\mu+1)\pi/4$$
for some $\mu\ge \nu$; in particular, the sector  $\bar\theta\le\arg z<(2\nu+1)\pi/4$ is `pole-free' for $w$. In a similar way one can show that the sector $(2\nu-1)\pi/4<\arg z\le \bar\theta$,
hence $\Sigma_\nu$ is `pole-free' for $w$.  The proof in the other cases runs along the same lines. \Ende

\subsection{\it Strings of poles} We will several times have to deal with sequences $(p_k)$ that tend to infinity and satisfy an approximative iterative scheme
\be{STRING}p_{k+1}=p_k+(\omega+o(1)) p_k^{-\tau},\ee
where $\omega\ne 0$ is  complex and $\tau=s/t>-1$ is  rational. Any such sequence is called a {\it string.}

\begin{lem}\label{lemma5}Any sequence ${\p}=(p_k)$ satisfying the recursion {\rm (\ref{STRING})} also satisfies
\begin{itemize}
\item[{\bf 1.}] $p_k=\big(k(1+\tau)\omega\big)^{t/(s+t)}(1+o(1))$ and
\item[{\bf 2.}] $(s+t)\arg p_k=t\arg\omega+o(1)$ $\mod 2\pi$, both times as $k\to\infty$, and
\item[{\bf 3.}] has counting function $n(r,{\p})=\ds\frac{r^{1+\tau}}{(1+\tau)|\omega|}(1+o(1))$  as $r\to\infty$.
\end{itemize}
\end{lem}

\beweis Writing $q_k=p_k^{1+\tau}$ it follows from (\ref{STRING}) that
$$q_{k+1}=q_k(1+(\omega+o(1))q_k^{-1})^{1+\tau}=q_k+(\tau+1)\omega+o(1),$$
hence {\bf 1.}\ follows from $q_k=(\tau+1+o(1))\omega k$, while assertions {\bf 2.}\ and {\bf 3.}\ are immediate consequences. \Ende

\begin{rem}\rm For $\tau$ some positive integer  the analogy with the dynamics of the rational map $R(p)=p+\omega p^{-\tau}$ is evident (see~\cite{NStDYN}
or any other text on rational dynamics); $R$ has a {\it parabolic fixed point} at infinity with $\tau+1$ {\it invariant
petals}, and the iterates $R^n(p)$ converge to infinity asymptotically to the rays $(\tau+1)\arg p=\arg\omega$ mod $2\pi$.\end{rem}

\begin{rem}\rm It follows from the proof of Proposition~\ref{RAYTOSECTOR} that `pole-free' sectors are `bordered' by strings of poles. In case
(\ref{ASYMPS})(ii)$_a$ and (iv)$_a$ the strings on both `sides' have the interlacing property: `between' any two poles with residue $1$ there is one with residue $-1$.
Otherwise the poles on both strings have one and the same residue.\end{rem}

\section{\bf Asymptotic Expansions}\label{ASYMPEXPANS}

\subsection{\it `Pole-free' sectors and asymptotic expansions} A function $f$ that is meromorphic on some sector $S$ is said to have an asymptotic expansion
$$f(z)\sim \sum\limits_{\nu=m}^\infty a_\nu z^{-\nu/q}\quad{\rm on~}S,$$
if $f(z)-\sum\limits_{\nu=m}^n a_\nu z^{-\nu/q}=\ko{-n/q}$
holds for  $n=m,m+1,\ldots$ as $z\to\infty$, uniformly on every proper sub-sector.
It is obvious that $S$ is `pole-free' if $f$ has an asymptotic expansion on $S$. For {\sc Painlev\'e} transcendents the converse is also true.

\begin{thm}\label{ASYMPEX}Every {\sc Painlev\'e} transcendent with `pole-free' sector $S$ has an asymp\-totic expansion on $S$.\end{thm}

\beweis Again we consider fourth transcendents only. Re-scaling along any sequence $(h_n)$ with $h_n\to\infty$ on any closed sub-sector $\tilde S\subset S$ yields limit functions without poles, hence constant limit functions either $\w=0$, $\w=-2$, or else $\w=-2/3$. This leads to the `leading terms' $\ko{}$, $-2z+\ko{}$, and $-2z/3+\ko{}$ of the asymptotic expansions in question.
We set $zv(z)=w(z)$ to obtain
$$z^{-2}(vv''-v'^2)=z^{-4}(2\beta+v^2)-4\alpha z^{-2} v^2+v^2(3v+2)(v+2)$$
with $c_0=\lim_{z\to\infty}v(z)\in\{0,-2,-2/3\}$. Theorem~\ref{AsympEx} in Appendix B applies if $c_0=-2$ and also if $c_0=-2/3$,
yielding $w\sim z\sum\limits_{k=0}^\infty c_kz^{-k}$ on $S$, but does not immediately apply if $c_0=0$.
In the latter case we note that $w=\ko{}$ implies $w'=o(1)$ and $w''=\ko{-1}$, thus $4z^2w^2$ dominates the terms $3w^4$, $8zw^3$, and
$-4\alpha w^2$ in equation (IV) (they are $o(|zw|^2$), and thus has to balance the term $-2\beta+2ww''-w'^2=-2\beta+o(1)$, that is,
$4z^2w^2+2\beta=o(|z|^2|w|^2)+o(1)$ holds and $zw$ tends to $\sqrt{-\beta/2}=\gamma$. In that case we set $v(z)=zw(z)$ to obtain
$$z^{-2}(vv''-v'^2)=z^{-2}(4\alpha v^2-8v^3)+z^{-4}(v^2-3v^4)+4v^2+2\beta.$$
Then for $\gamma=\sqrt{-\beta/2}\ne 0$ again Theorem~\ref{AsympEx} applies with $c_0=\lim_{z\to\infty}v(z)=\gamma$.
Finally, for $\gamma=0$, we have $w=\ko{-1}$. If $w=\ko{-n}$ is assumed to be true for some $n\ge 1$ we obtain
$$\begin{array}{rcl}4z^2w^2&=&2ww''-w'^2-3w^4-8zw^3+4\alpha w^2\cr
&=&\ko{-2n-2}+\ko{-4n}+\ko{-3n+1}+\ko{-2n}=\ko{-2n},\end{array}$$
hence $w=\ko{-n-1}$ and $w\sim 0$ (all coefficients vanish). \Ende

\begin{rem}\rm In case (I) and (II) we set $z=t^2$ and $tv(t)=w(z)$ to obtain
$$\begin{array}{lrcl}
{\rm (i)}&\ds t^{-1}\ddot v+t^{-2}\dot v&=&\ds t^{-3} v+4(1+6v^2)\quad{\rm and}\cr
{\rm (ii)}&\ds t^{-1}\ddot v+t^{-2}\dot v&=&\ds 4\alpha t^{-3}+t^{-5}v+4v(1+2v^2),\end{array}$$
respectively ($~\dot{}$ means differentiation with respect to $t$), with
(i) $c_0=\lim_{t\to\infty}v(t)=\sqrt{-1/6}$ on $S^{1/2}$ and (ii) $c_0=\sqrt{-1/2}$ or $c_0=0$; in case of $c_0=0$
we use $w=\ko{1/2}$ and $w^3=o(|zw|)$ to obtain $zw+\alpha=o(1)$.\end{rem}

\subsection{\it Explicit expansions}It is not hard to confirm the following detailed expansions:
\be{asympexw}\begin{array}{lrcl}
{\rm (i)}_{\phantom{a}}&w(z)&=&\ds \sqrt{-z/6}-\frac1{48z^2}+\GO{-9/2}\cr
{\rm(ii)}_a&w(z)&=&\ds\phantom{-2z}-\frac{\alpha}z+\frac{2\alpha(\alpha^2-1)}{z^4}+\GO{-7}\cr
{\rm(ii)}_b&w(z)&=&\ds\sqrt{-z/2}+\frac\alpha{2z}+\GO{-5/2}\cr
{\rm(iv)}_a&w(z)&=&\ds-\frac23 z+\frac\alpha
 z-\frac{3\alpha^2-9\gamma^2+1}{4z^3}+\GO{-5}\cr
{\rm(iv)}_b&\ds w(z)&=&\ds-2z-\frac\alpha z+\frac{3\alpha^2-\gamma^2+1}{4z^3}+\GO{-5}\cr
{\rm(iv)}_{c}^\pm&w(z)&=&\ds\phantom{-2z}\pm\frac\gamma z-\frac{2\gamma^2\mp\alpha\gamma}{2z^3}+\GO{-5}
\qquad(\gamma^2=-\beta/2).\end{array}\ee
\subsection{\it The case $\beta=0$}The expansions (iv)$_{c}^\pm$ are only significant if $\gamma=\sqrt{-\beta/2}\ne 0$. We thus suppose $\beta=0$ and $w\sim 0$ on some sector $S$, which means
$w(z)=\ko{-n}$ as $z\to\infty$ on $S$ for every $n\in\N$. The logarithmic derivative $y=w'/w$ satisfies
$$y'=P(z)-\frac12 y^2\quad{\rm with}~P(z)=2z^2-2\alpha+\frac32zw(z)+2w(z)^2\sim 2z^2-2\alpha.$$
From the considerations in \cite{NStRiccati} it then follows that $y=w'/w$ has an asymptotic expansion $y\sim \pm 2z+\cdots$ on $S$, and $w\to 0$ on $S$ requirers $\Re(\pm z^2)<0$ on $S$
(for some sign), that is,
\be{asympexw0}\begin{array}{rrcl}
{\rm (iv)}_c^-&\ds\frac{w'}w&\sim&\ds -2z+\frac{\alpha-1}z+\frac{\alpha^2-4\alpha+3}{4z^3}+\cdots\quad{\rm if~}S\subset\Sigma_0\cup\Sigma_2\cr
{\rm (iv)}_c^+&\ds\frac{w'}w&\sim&\ds \phantom{-}2z-\frac{\alpha+1}z-\frac{\alpha^2+4\alpha+3}{4z^3}+\cdots\quad{\rm if~}S\subset\Sigma_1\cup\Sigma_3.\end{array}\ee
A similar result holds in (\ref{asympexw})(ii)$_a$ if $\alpha=0$:
$\ds\frac{w'}w=\pm\sqrt z-\frac1{4z}\mp\frac 5{32z^2\sqrt z}+\cdots$ holds on sectors $\Re(\pm z\sqrt z)<0$.

\subsection{\it Asymptotics of the first integrals}From (\ref{FirstIntegral}) and (\ref{asympexw}) it easily follows that
\be{asympexW} \begin{array}{lrcl}
{\rm (i)}&W(z)&=&\ds\frac{1}9\sqrt{-6}z^{3/2}+\frac1{48z}+\GO{-7/2}\cr
{\rm(ii)}_a&W(z)&=&\left\{\begin{array}{lr}\ds-\frac{\alpha^2}z+\frac{\alpha^2(\alpha^2-1)}{z^4}+\GO{-7}&(\alpha\ne 0)\cr
\GO{M}e^{-\frac23|\Re z^{\frac 32}|}&(\alpha=0)\end{array}\right.\cr
{\rm(ii)}_b&W(z)&=&\ds -\frac14 z^2+\sqrt{-2}\alpha z^{1/2}+\frac{1+4\alpha^2}{8z}+\GO{-5/2}\cr
{\rm(iv)}_a&W(z)&=&\ds -\frac8{27}z^3+ \frac23\alpha z-\frac{3\alpha^2+9\gamma^2-1}{6z}+\GO{-3}\cr
{\rm(iv)}_b&\ds W(z)&=&\ds 2\alpha z+\frac{\alpha^2-\gamma^2+1}{2z}+\GO{-3}\cr
{\rm(iv)}_c^\pm&W(z)&=&\left\{\begin{array}{lr}\ds\pm 2\gamma z+\frac{\gamma^2\mp\alpha\gamma}{z}+\GO{-3}&(\gamma\ne 0)\phantom{.}\cr
\GO{M}e^{-|\Re z^2|}&(\gamma=0).\end{array}\right.
\end{array}\ee

\subsection{\it The Hastings-{McLeod} solution}Equation $w''=zw+2w^3$ has a unique solution, named after {\sc Hastings} and {\sc McLeod} \cite{HastingsMcLeod},
that decreases on the real line. Moreover, it satisfies $w(x)\sim 0$ as $x\to+\infty$ and $w(x)\sim\sqrt{-x/2}$ as $x\to-\infty$. From
Proposition~\ref{RAYTOSECTOR} and Theorem~\ref{ASYMPEX} it follows that the asymptotic expansions $w\sim 0$ and $w(z)\sim\sqrt{-z/2}+\cdots$ hold throughout $|\arg z|<\pi/3$ and $|\arg z-\pi|<\pi/3$,
respectively. The poles are asymptotically restricted to the sectors $|\arg z\mp\pi/2|<\pi/6$. Writing $w''=(z+2w(z)^2)w$ and noting that $w(z)=\ko{-n}$  on $|\arg z|<\pi/3$
for every $n\in\N$ yields $w(z)\sim k{\sf Ai}(z)$ for some real constant $k$, actually $k=1$; {\sf Ai} denotes the {\sc Airy} function.
More general, for $0<k<1$ there exists a unique solution, named after
{\sc Ablowitz} and {\sc Segur} \cite{AbloSegur1, AbloSegur2}, that is bounded on
$\R$ and satisfies $w(x)\sim k{\sf Ai}(x)$ as $x\to +\infty$. Since in that case $w$ has the asymptotics $w(x)=(-x)^{-1/4}$ times an oscillating term as $x\to -\infty$,
$w$ cannot have an asymptotic expansion on some sector about $\arg z=\pi$, and hence the poles of $w$ must accumulate at $\arg z=\pi$.
Also Proposition~\ref{RAYTOSECTOR} is void in the present case since $\arg z=\bar\theta=\pi$ is a {\sc Stokes} ray.

\subsection{\it The {Clarkson}-{McLeod} solution}For $\beta=0$ and $\alpha$ real it is conjectured that there exists a unique real solution to equation (IV)---the so-called
{\it {\sc Clarkson}-{\sc McLeod} solution} \cite{ClarkMcLeod,ItsKapaev}---satisfying
$w(x)\sim 0$ as $x\to+\infty$ and $w(x)\sim -2 x+\cdots$ as $x\to -\infty$. If it exists, and existence is supported by numerical experiments and by analogy to the
{{\sc Hastings}-{\sc McLeod} solution}, then
the corresponding asymptotic expansions $w\sim 0$ (even $w'/w\sim-2z+\cdots$ by (\ref{asympexw0})) and $w\sim -2z+\cdots$
hold on $|\arg z|<\pi/4$  and $|\arg z-\pi|<\pi/4$, respectively.

\subsection{\it B\"acklund transfomations}\label{BAECKLT}
Generally speaking, a {\it {\sc B\"acklund} transformation} is any change
of variables that transforms some given differential equation into itself or into
a differential equation of the same type.  We confine ourselves to equation (IV). The {\sc B\"acklund} transformation due to {\sc Lukashevich} \cite{Luka} (see also \cite{GLS})
 \be{BackIV}\tilde w=\frac{w'-2\gamma-2zw-w^2}{2w}\quad(\gamma^2=-\beta/2)\ee
transforms equation (IV) into the same equation with new parameters
$$\tilde \alpha=(1-\alpha+3\gamma)/2\quad{\rm and}\quad\tilde \gamma=(1+\alpha+\gamma)/2;$$
the (pointwise) inverse transformation is given by
\be{BackIVinverse}\ds w=-\frac{\tilde w'+2\tilde \gamma+2z\tilde w+\tilde w^2}{2\tilde w}\quad(\tilde \gamma^2=-\tilde \beta/2)\ee
with $\alpha=(-1-\tilde \alpha+3\tilde \gamma)/2$ and $\gamma=(-1+\tilde \alpha+\tilde \gamma)/2$; we note that
\be{ALPHAGAMMA}\tilde \alpha-\tilde \gamma=-(\alpha-\gamma)\ee
holds in any case.
Applying (\ref{BackIV}) twice yields equation (IV) with parameters $\alpha+1$ and $\gamma+1$.
{\sc Lukashevich}'s transformations form the master pair (actually, quadruple: $\gamma$ and $\tilde \gamma$ may be replaced by $-\gamma$ and $-\tilde \gamma=(1+\alpha-\gamma)/2$, respectively), since every (known and unknown) {\sc B\"acklund} transformation may be expressed by certain compositions of
(\ref{BackIV}) and (\ref{BackIVinverse}), see \cite{BCH}.

\subsection{\it Invariance}We close this section by proving the invariance of cluster sets, `pole-free' sectors, asymptotic expansions and
the {\sc Nevanlinna} characteristic under {\sc B\"acklund} transformations.

\begin{prp}\label{BACKCHANGE}For every fourth transcendent the cluster set, `pole-free' sectors, and the {\sc Nevanlinna} characteristic remain
invariant under the corresponding {\sc B\"acklund} transformations {\rm (\ref{BackIV})} and  {\rm (\ref{BackIVinverse})},
the latter meaning $T(r,\tilde w)\asymp T(r,w)$. Moreover, the {\sc B\"acklund} transformation
{\rm (\ref{BackIV})} changes
$$w\sim\left\{\begin{array}{l} -\frac23z+\cdots\cr-2z+\cdots\cr\end{array}\right.~{\rm into}~
\tilde w\sim\left\{\begin{array}{l}-\frac23z+\cdots\cr \tilde\gamma/z+\cdots\end{array}\right.$$
For $\gamma\ne 0$ it also changes
$$w\sim\left\{\begin{array}{l} -\gamma/z+\cdots\cr \gamma/z+\cdots\cr\end{array}\right.~{\rm into}~
\tilde w\sim\left\{\begin{array}{l}-\tilde\gamma/z+\cdots\cr -2z+\cdots\end{array}\right.$$
while for $\gamma=0$, $w\sim 0$ $(w\not\equiv 0)$ is changed into
 $$\tilde w\sim\left\{\begin{array}{ll} -2z+\cdots&{\rm on~} \Sigma_0 {\rm ~and~}\Sigma_2,\cr
-\tilde\gamma/z+\cdots&{\rm on~} \Sigma_1 {\rm ~and~}\Sigma_3.\end{array}\right.$$
\end{prp}

\begin{rem}\rm A similar result holds for second transcendents.\end{rem}

\beweis The assertion about the cluster set follows from
$$w'-2\gamma-2zw-w^2=2w\tilde w=-(\tilde w'+2\tilde\gamma+2z\tilde w+\tilde w^2),$$
hence $w'+\tilde w'-2(\gamma-\tilde\gamma)=2zw+w^2-(2z\tilde w+\tilde w^2)=W'-\tilde W'$ and
\be{wundW}W-\tilde W=w+\tilde w-2(\gamma-\tilde\gamma)z+const.=\GO{}\quad(z\notin\P_\delta(w)\cup\P_\delta(\tilde w)).\ee
Moreover, from $\ds \tilde w=\frac{w'-2\gamma}{2w}-z-\frac w2$ it follows that
$$T(r,\tilde w)+O(\log r)=N(r,\tilde w)\le N(r,1/w)+N(r,w)\le 2T(r,w)+O(1),$$
and, in the same manner, $T(r,w)\le 2T(r,\tilde w)+O(\log r)$.
To prove the statements about the change of the asymptotic expansions is just a matter of elementary computations, at least if $\beta\ne 0$.
For $\gamma=0$ it is also true that the expansion (iv)$_b$ is transformed into (iv)$_c^+$. The other cases need a more subtle argument. From (\ref{asympexw0}), $\gamma=0$,
$w\sim 0$, and $\ds \tilde w=\frac{w'}{2w}-z-\frac w2$ we obtain $\tilde w\sim
\ds-2z+\frac{\alpha-1}{2z}+\cdots$ on $\Sigma_0\cup\Sigma_2$,
and $\tilde w\sim \ds-\frac{\alpha+1}{2z}+\cdots$ on $\Sigma_1\cup\Sigma_3.$ \Ende

\section{\bf Truncated Painlev\'e Transcendents}\label{MULTTRUNC}

\subsection{\it Truncated solutions}{\sc Painlev\'e} transcendents having some `pole-free' sector $S:|\arg z-\theta_\nu|<\epsilon$ about the {\sc Stokes} ray $\sigma_\nu:\arg z=\theta_\nu$
are called {\it truncated} (along $\sigma_\nu)$.
By Theorem~\ref{ASYMPEX}, $w$ has an asymptotic expansion on $S$, and from Proposition~\ref{RAYTOSECTOR} applied to $\theta_\nu-\epsilon<\bar\theta<\theta_\nu$ and
$\theta_\nu<\bar\theta<\theta_\nu+\epsilon$ it follows that the asymptotic expansion extends to $\Sigma_\nu\cup\sigma_\nu\cup\Sigma_{\nu+1}$. In particular, any  transcendent $w$ that is truncated along $\sigma_\nu$ has the same asymptotic expansion on adjacent {\sc Stokes} sectors
$\Sigma_\nu$ and $\Sigma_{\nu+1}$.(\footnote{In cases (i) and (ii) ``the same expansion'' means that the square-root $\sqrt{-z}$ is assumed to be continuous across $\sigma_\nu.$}) The converse is also true:

\begin{prp}\label{ADJACENT}Let $w$ be any {\sc Painlev\'e} transcendent that has the same asymptotic expansion {\rm (\ref{asympexw})} on adjacent {\sc Stokes} sectors
$\Sigma_\nu$ and $\Sigma_{\nu+1}$. Then $w$  is truncated along $\sigma_\nu$, and the asymptotic expansion
holds on $\Sigma_\nu\cup\sigma_\nu\cup\Sigma_{\nu+1}$.\end{prp}

\beweis In any case $\ds F(z)=e^{-\int W(z)\,dz}$
is an entire function of finite order with simple zeros at the poles of $w$. Since again the idea of proof is the same in all cases we will consider this time
(\ref{asympexw})(iv)$_b$, say,
with $\Sigma_\nu:|\arg z-\nu\pi/2|<\pi/4$ and
$$W(z)=\ds 2\alpha z+\frac{\alpha^2-\gamma^2+1}{2z}+\GO{-3},$$
which holds by (\ref{asympexW}), hence
$$H(z)=F(z)e^{\alpha z^2}z^{(\alpha^2-\gamma^2+1)/2}=
\left\{\begin{array}{rl}c_\nu +o(1)&{\rm as~}z\to\infty~{\rm on~} \Sigma_\nu\cr
c_{\nu+1}+o(1)&{\rm as~}z\to\infty~{\rm on~} \Sigma_{\nu+1}\end{array}\right.$$
holds for some complex constants $c_\nu$ and $c_{\nu+1}$. We set
$$f(\zeta)=H(z)\quad{\rm with}~ z=e^{i(2\nu+1)\pi/4}\sqrt{\zeta}\hspace{3mm}{\rm and~}~\sqrt\zeta>0 {\rm~on~} \zeta>0.$$
Since $f$ has finite order $\ds\limsup_{\zeta\to\infty}\frac{\log\log|f(\zeta)|}{\log|\zeta|}\le 2$ on $|\arg \zeta|<\pi-\eta$ and limits $c_\pm$ as $\zeta=re^{\pm i\delta}\to\infty$
for every $0<\delta<\pi$, the {\sc Phragm\'en}-{\sc Lindel\"of} Principle (see Theorem~\ref{PhragLind} in Appendix~A) yields $c_-=c_+=c$ and $f(\zeta)\to c$ on
$|\arg\zeta|\le\delta$. From Proposition~\ref{RAYTOSECTOR} it then follows that the half-plane $|\arg z-(2\nu+1)\pi/4|<\pi/2$ is `pole-free' for $w$, and $w$ has the asymptotic expansion  (\ref{asympexw})(iv)$_b$ there. \Ende

\subsection{\it Existence and uniqueness of truncated solutions}We will first prove a result valid in all cases, which will afterwards be refined for equation (I).

\begin{thm}\label{EXUNIQUE}{\rm (\cite{Novok1, Novok2} and many others)} Given any appropriate {\sc Stokes} ray $\sigma_\nu$, there exist
\begin{itemize}\item[(i)] first transcendents that are truncated along $\sigma_\nu$ and $\sigma_{\nu+1}$,
\item[(ii)] solutions to (II) that are truncated along $\sigma_\nu$, and
\item[(iv)]  solutions to (IV) that are truncated along $\sigma_\nu$;\end{itemize}
in {\rm (ii)} and {\rm (iv)} rational solutions are included. The asymptotic expansion {\rm (\ref{asympexw})} along $\sigma_\nu$, hence on $\Sigma_\nu\cup\sigma_\nu\cup\Sigma_{\nu+1}$ may be prescribed; it uniquely determines the solutions in question.\end{thm}

\beweis Again we restrict ourselves to fourth transcendents. By Theorem~\ref{EXISTENCE} in Appendix~A there exists some fourth transcendent $w$ with prescribed asymptotics (iv) on any given sector $S$ of central angle $\pi/2$. In particular, if $S$ contains some {\sc Stokes} ray $\sigma_\nu$, the asymptotic expansion holds on  $\Sigma_\nu\cup\sigma_\nu\cup\Sigma_{\nu+1}$ by Proposition~\ref{RAYTOSECTOR},
and $w$ is truncated along $\sigma_\nu$.  The proof is similar in the other cases on combination of  Theorem~\ref{EXISTENCE} and Proposition~\ref{RAYTOSECTOR}.

To prove uniqueness we consider fourth {\sc Painlev\'e} transcendents $w_1$ and $w_2$  having the same asymptotic expansion (\ref{asympexw}) on some
sector $S$. Then $u=w_1-w_2$ tends to zero as $z\to\infty$ more than geometrically fast on every proper sub-sector of $S$, that is, for every positive integer $n$, $z^nu$ tends to zero. In all cases
we will derive a linear differential equation
\be{LINI}y''+Q(z)y=0,\ee
which on one hand has no non-trivial solution that tends to zero as $z\to\infty$ on $S$ if the sector $S$ has central angle $\Theta>\pi/2$, and on the other hand has some non-trivial solution tending to zero provided $u=w_1-w_2\not\equiv 0$. This will show that $w_1=w_2$.

From
$$\begin{array}{rcl}w_1w_1''-w_2w_2''&=&w_1u''+w_2''u=w_2u''+w_1''u, {\rm~hence}\cr
2w_1w_1''-2w_2w_2''&=&(w_1+w_2)u''+(w_1''+w_2'')u\cr
w_1'^2-w_2'^2&=&(w_1'+w_2')u'\cr
w_1^4-w_2^4&=&(w_1^3+w_1^2w_2+w_1w_2^2+w_2^3)u~{\rm ~etc.}\end{array}$$
we obtain
$(w_1+w_2)u''+(w_1''+w_2'')u=(w_1'+w_2')u'+h(z)u$
with
$$h(z)=3(w_1^3+w_1^2w_2+w_1w_2^2+w_2^3)+8z(w_1^2+w_1w_2+w_2^2)+4(z^2-\alpha)(w_1+w_2),$$
hence the linear differential equation
\be{LINIV}u''+a(z)u'+b(z)u=0\ee
with coefficients~ $a(z)=\ds-\frac{w_1'(z)+w_2'(z)}{w_1(z)+w_2(z)}$ and~ $b(z)=\ds\frac{w_1''(z)+w_2''(z)-h(z)}{w_1(z)+w_2(z)}.$
In case of (\ref{asympexw})(iv)$_c$ we first assume in addition that $\gamma\ne 0$ (otherwise we would have to divide by $w_1+w_2\sim 0$ which would cause difficulties).
The substitution
$$y=u\,e^{\frac12\int a(z)\,dz}=(w_1(z)+w_2(z))^{-1/2}u=z^{\epsilon/2}(c+o(1))u$$
($c\ne 0$, $\epsilon=1$ in case (iv)$_c^\pm$, and $\epsilon=-1$ otherwise)
transforms (\ref{LINIV}) into (\ref{LINI}) with
$Q(z)= b(z)-\frac14 a(z)^2-\frac12 a'(z)=b(z)+\GO{-2},$
and, in more detail,
$$Q(z)=\left\{\begin{array}{ll}\phantom{+}\frac43z^2+\GO{-2}&{\rm (iv)}_a\cr
-4z^2-8\alpha+\GO{-2}&{\rm (iv)}_b\cr -4z^2+4\alpha\mp 12\gamma+\GO{-2}&{\rm (iv)}_c^\pm\end{array}\right.$$
The asymptotic integration of (\ref{LINI}) is well understood (see {\sc Hille}~\cite{EH} or else {\sc Wasow}~\cite{Wasow1}). We just need a qualitative result:
writing $Q(z)=cz^\tau+\ko{\tau}$, every non-trivial
solution to (\ref{LINI}) tends to infinity exponentially on some sub-sector of $S$ if the central angle $\Theta$ of $S$ is greater than
$\frac{2\pi}{\tau+2}$. In our case this yields $y\equiv 0$ if $\Theta>\pi/2$, and this proves
Theorem~\ref{EXUNIQUE} up to the case (iv)$_c^\pm$ with $\gamma=0$. Here we compare the solution $w_1=w\sim 0$ with the trivial solution $w_2\equiv 0$.
From $w(z)=\ko{-n}$ for every $n\in\N$, and $(w'(z)/w(z))^2=4z^2-4\alpha\pm 4+\GO{-2}$, which follows from (\ref{asympexw0}), we obtain
$$-Q(z)=\frac{w'(z)^2}{2w(z)^2}+\frac 32 w(z)^2+2zw(z)+2z^2-2\alpha=4z^2-4\alpha\pm 2+\GO{-2},$$
hence essentially the same equation as in case of (iv)$_c^\pm$ with $\gamma\ne 0$.
The proof of uniqueness is much easier for equations (I) and (II); we obtain $u''+Q(z)u=0$ with $\tau=1/2$ and $\tau=1$, hence either $u\equiv 0$
or else $\Theta\le \frac 45\pi$ and $\Theta\le \frac23\pi$, respectively. \Ende

\subsection{\it Triply truncated first transcendents}
 Following {\sc Boutroux}~\cite{Boutroux1},  {\sc Painlev\'e}'s equation (I) has five {\it triply truncated} solutions
(also called {\it tritronqu\'ee}); for a recent existence proof see {\sc Joshi} and {\sc Kitaev} \cite{JosKia}. Since (I) is invariant under the transformation $w\mapsto a^2w(az)$ with $a^5=1$, it
suffices to prove the existence of a triply truncated solution $w_0$ having the asymptotics
$w_0(z)\sim -\sqrt{-z/6}$ with $\Re\sqrt{-z/6}>0$ on $|\arg z-\pi|<4\pi/5$, which thus is truncated along $\sigma_1$, $\sigma_2$, and $\sigma_3$.
By Theorem~\ref{EXUNIQUE}(i) there exist uniquely determined first transcendents $w_1$ and $w_2$ with asymptotics
$w_{1,2}= -\sqrt{-z/6}+\GO{-2}$ and $\Re\sqrt{-z/6}>0$
on $\pi/5<\arg z<7\pi/5$ and on $3\pi/5<\arg z<9\pi/5$, respectively. Then again $y=w_1-w_2$ satisfies
$$y''=6(w_1(z)+w_2(z))y=(-12\sqrt{-z/6}+\GO{-2})y$$
on $3\pi/5<\arg z<7\pi/5.$ We set $z=-6t^{4/5}$ and
$y(z)=t^{-1/10}v(t)$ on the right half-plane $\Re t>0$ ($t^{4/5}>0$, $\sqrt{-z/6}=t^{2/5}>0$ and $t^{-1/10}>0$ if $t>0$) to obtain
\be{DGLv}{100}\ddot v+\big[27648+9t^{-2}+O(|t|^{-11/5})\big]v=0.\ee
Now every non-trivial solution to (\ref{DGLv}) tends to infinity exponentially as $t\to\infty$
at least on one of the sectors $\delta<\arg t<\pi/2-\delta$ and  $-\pi/2+\delta<\arg t<-\delta$.
This proves $v\equiv 0$ and $y\equiv 0$, and $w_0=w_1=w_2$ is truncated along $\sigma_1$, $\sigma_2$, and $\sigma_3$.

\begin{rem}\rm We note that $\overline{w_0(\bar z)}$ is also a solution that is truncated along $\sigma_1$, $\sigma_2$, and $\sigma_3$, and hence coincides with $w_0$ by uniqueness. In other words, the tritronqu\'ee solution is real on the real axis.
We also note that if $w$ is any solution that is truncated along $\sigma_1$ and $\sigma_2$ resp.\  $\sigma_2$ and $\sigma_3$ coincides with $w_0$ by Theorem~\ref{EXUNIQUE}
if it has the same asymptotics as $w_0$ on $\Sigma_1\cup\sigma_1\cup\Sigma_2$ resp.\ $\Sigma_2\cup\sigma_2\cup\Sigma_3$.\end{rem}

\begin{rem}\rm In contrast to Theorem~\ref{EXUNIQUE} the triply truncated solution $w_0$ is not only unique but also determines its asymptotics. In other words, there is no solution $\tilde w_0$ that is truncated along $\sigma_1$, $\sigma_2$, and $\sigma_3$ and is asymptotic to $+\sqrt{-z/6}$ with $\Re\sqrt{-z/6}>0$.
Otherwise we will consider the solution $v(z)=e^{4\pi i/5}\tilde w_0(e^{2\pi i/5}z)$, which is truncated along
the {\sc Stokes} rays $\sigma_0$, $\sigma_1$, and $\sigma_2$. Then $w_0$, $\tilde w_0$ and $v$ have asymptotic expansions on $\pi/5<\arg z<7\pi/5$,
which by Theorem~\ref{EXUNIQUE} is only possible if two of these solutions agree, which obviously is not the case.\end{rem}

\section{\bf The Dichotomy of the Order}\label{DICHOTOMY}
On combination with work of {\sc Mues} and {\sc Redheffer}~\cite{MuesRedheffer1} and {\sc Shimomura}~\cite{Shimomura0}
it follows that first transcendents have order of growth $\frac52$.  The situation is quite different for second and fourth transcendents, they have order
of growth between $\frac 32$ and $3$ and between $2$ and $4$, respectively.
More precisely they satisfy
$$C_1r^{\frac32}\le T(r,w)\le C_2r^3\quad{\rm and}\quad C_1r^{2}\le T(r,w)\le C_2r^4,$$
in the respective cases, see \cite{HL2,HL4,Shimomura0,NStBoutroux}.

\subsection{\it Weber-Hermite solutions}\label{WeberHermite}The solutions to the {\it {\sc Weber-Hermite} equations}
\be{WEBERHERMITE}w'=2\gamma\pm(2zw+w^2)\quad(\gamma=-1\mp\alpha)\ee
also solve equation (IV) with $\beta=-2(1\pm\alpha)^2$, and satisfy $n(r,w)\asymp r^2$.(\footnote{More precisely, $\ds n(r,w)=\nu(w)\frac{r^2}{2\pi}+o(r^2)$ holds, were either $\nu(w)=4$ or else $\nu(w)=2$ (\cite{NStRiccati}).}) Repeated application of the {\sc B\"acklund}
transformations (\ref{BackIV}), (\ref{BackIVinverse}), and the transformations $w(z)\mapsto-iw(iz)$ and $w\mapsto\overline{w(\bar z)}$ to the solutions to
(\ref{WEBERHERMITE}) leads to the so-called {\it {\sc Weber-Hermite} solutions}. Just like the solutions to (\ref{WEBERHERMITE})
they have counting function of poles $n(r,w)\asymp r^2$.
In section~\ref{SUBNORMALSOLUTIONS} it will be shown that the {\sc Weber-Hermite} solutions coincide with those satisfying $n(r,w)=O(r^2)$.
In case of (II) the role of the {\sc Weber-Hermite} equations is taken by the {\it {\sc Airy} equations}
$w'=\pm(z/2+w^2)$ with $n(r,w)\asymp r^{3/2}.$

\subsection{\it Solutions of maximal order of growth}It is not hard to identify the second and fourth transcendents of order $3$ and $4$ mean type, respectively.

\begin{thm}\label{MaximalOrder}Let $w$ be any {\sc Painlev\'e} transcendent whose cluster set contains some non-exceptional parameter $\c_0$. Then $w$ has maximal order
of growth. More precisely, there exist sequences $r_k\to\infty$  such that
$${\rm (i)}\quad n(r_k,w)\asymp r_k^{5/2},\quad{\rm (ii)}\quad n(r_k,w)\asymp r_k^{3},\quad{\rm and\quad (iv)}\quad n(r_k,w)\asymp r_k^{4}$$
holds in the respective case.\end{thm}

\begin{rem}\rm {\sc Hinkkanen} and {\sc Laine} \cite{HL3} for equation (II), and {\sc Sriponpaew} \cite{BS} and {\sc Classen}~\cite{Classen} for equation (IV) provided different proofs of Theorem~\ref{MaximalOrder}.\end{rem}

\beweis Again the idea is the same in all cases, as usual we prefer to consider fourth transcendents. There exists some sequence $(h_k)$ such that
$h_k^{-3}W(h_k)$ tends to some constant $\c_0\ne 0, -8/27$. Since $z^{-3}W(z)$ varies slowly, this following from the fact that
$$\frac d{dz} z^{-3}W(z)=-3z^{-4}W(z)+z^{-3}(2zw(z)+w(z)^2)=\GO{-1}$$
holds outside $\P_\delta$, there exists some $\epsilon>0$ such that
$$|z^{-3}W(z)-\c_0|<\ts\frac12\min\{|\c_0|,|\c_0+8/27|\}$$
holds on $D_k=\{z:|z-h_k|<\epsilon|h_k|\}\setminus\P_\delta;$ $\delta>0$ is fixed. In other words, for every appropriate sequence $(\tilde h_\nu)$ with $\tilde h_\nu\in D_{k_\nu}$ the limit function $\w=\lim_{\nu\to\infty}w_{\tilde h_\nu}$ is doubly periodic and solves some differential equation
$\w'^2=\w^4+4\w^3+4\w^2-4\c\w$ with parameter satisfying $|\c-\c_0|\le\frac12\min\{|\c_0|,|\c_0+8/27|\}$, hence $|\c|\ge \frac12|\c_0|$ and also $|\c+8/27|\ge \frac12|\c_0+8/27|$.
The corresponding lattice $\Lambda_\c$ has a fundamental parallelogram whose diameter is bounded independent of $\c$, hence there exist $R>0$ and $k_0\in\N$,
such that for every $k\ge k_0$ and center $z_0$ satisfying $|z_0-h_k|<\epsilon|h_k|$, the disc ${\sf D}_R(z_0)=\{z:|z-z_0|<R|h_k|^{-1}\}$  contains at least one pole of $w$. Now the disc $|z-h_k|<\epsilon|h_k|$ contains $\asymp |h_k|^4$ centers $z_\nu$ of mutually disjoint discs ${\sf D}_R(z_\nu)$ (it is the same to say that, as $r\to\infty$, the disc $|\z|<r$ contains $\asymp r^2$ integers $m+in$: just map $|\z|<\frac\epsilon{2R}|h_k|^2$ onto $|z-h_k|<\epsilon|h_k|$ by $\z\mapsto h_k+2R|h_k|^{-1}\,\z$), and since the discs ${\sf D}_R(z_\nu)$ are contained in $|z-h_k|<2|h_k|$, say, it follows that
$n(2|h_k|,w)\ge \mu|h_k|^4$
holds for some $\mu>0$ that is independent of $k$. Together with $n(r,w)=O(r^4)$ this yields $n(r_k,w)\asymp r_k^4$ for the sequence $r_k=2|h_k|\to\infty$. \Ende

\begin{rem}\label{StrongConj}\rm There is some evidence to believe that a much stronger result holds: {\it Given any fourth [second] transcendent $w$ and any {\sc Stokes} sector $\Sigma$, then
either $w$ has an asymptotic expansion on $\Sigma$ or else the cluster set $\CL_{\Sigma^*}(w)$ of $w$ restricted to any closed sub-sector $\Sigma^*$ of $\Sigma$
contains neither $0$ nor $-8/27$ [neither $0$ nor $-1/4$].}
In other words, {\it either each or else none limit function $\w=\lim_{h_n\to\infty}w_{h_n}$ with $(h_n)\subset\Sigma^*$ is constant.}\end{rem}

\section{\bf Sub-normal Solutions}\label{SUBNORMALSOLUTIONS}
Fourth transcendents satisfying $n(r_k,w)=O(r_k^2)$ on some sequence $r_k\to\infty$ are called {\it sub-normal}. It will soon turn out that
sub-normal solutions even satisfy $n(r,w)=O(r^2)$ as $r\to\infty$.  We also note that by Theorem~\ref{MaximalOrder}, fourth transcendents
satisfying $n(r,w)=o(r^4)$ have cluster set either $\CL(w)=\{0\}$ or else $\CL(w)=\{-8/27\}$. In particular, this is true for sub-normal solutions.

\subsection{\it Subnormal solutions with cluster set $\CL(w)=\{-8/27\}$} Similar to the case of equation (II), which was considered in \cite{NStSubnormal}, it turns out that the following is true.

\begin{thm}\label{SUBNORMALIVnot}Sub-normal fourth {\sc Painlev\'e} transcendents  with cluster set $\CL(w)=\{-8/27\}$ do not exist.\end{thm}

\beweis Let $w$ be any fourth transcendent with  counting function of poles $n(r_k,w)=O(r_k^2)$ for some sequence $r_k\to\infty$, and cluster set  $\CL(w)=\{-8/27\}$. The re-scaling method then shows that every
pole of sufficiently large modulus is contained in some uniquely determined string of poles $\p=(p_k)$ satisfying
$$p_{k+1}=p_k\pm(\pi\sqrt 3+o(1))p_k^{-1},$$
hence $p_k$ is given by
$p_k=(\pm(2\pi\sqrt 3+o(1))k)^{\frac 12}$ with fixed sign $\pm$ and fixed square-root (see Lemma~\ref{lemma5}). From our hypothesis
and
$$\ds n(r,\p)=\frac{r^2}{2\pi\sqrt 3}+o(r^2)$$
it follows that there are only finitely many
such strings, each of them being asymptotic to some ray $\arg z=\nu\pi/2$. Thus the corresponding {\sc Stokes} sectors
$\Sigma_\nu$ are `pole-free', and $w$ has one and the same asymptotic expansion (\ref{asympexw})(iv)$_a$ on each $\Sigma_\nu$. By Proposition~\ref{ADJACENT}
the asymptotic expansion (\ref{asympexw})(iv)$_a$ holds on the whole plane, hence $w$ is a rational function. This proves Theorem~\ref{SUBNORMALIVnot}. \Ende

\subsection{\it Subnormal solutions with cluster set $\CL(w)=\{0\}$}We shall now prove the main theorem on sub-normal solutions. The corresponding result for second transcendents has been proved in \cite{NStSubnormal}; Theorems~\ref{SUBNORMALIVnot} and \ref{SUBNORMALIV} have also been proved by {\sc Cla{ss}en}~\cite{Classen}.

\begin{thm}\label{SUBNORMALIV}The sub-normal fourth transcendents with cluster set $\CL(w)=\{0\}$ coincide with the {\sc Weber-Hermite} solutions.\end{thm}

\beweis  Let $w$ be any sub-normal fourth transcendent with cluster set $\CL(w)=\{0\}$. Again the re-scaling method,
on combination with Lemma~\ref{lemma5} with $\tau=1$ and $\omega=\pm \pi i$ shows that the set of poles of $w$ on $|z|>r_0$ consists of finitely many
strings $\p=(p_k)$, each being asymptotic to some ray $\arg z=(2\nu+1)\pi/4$ $(0\le\nu\le 3)$ and having counting function $\ds n(r,\p)=r^2/2\pi+o(r^2)$.
This implies
\be{SUBNORMAL}n(r,w)=\nu(w)\frac{r^2}{2\pi}+o(r^2),\ee
where $\nu(w)$ denotes the number of strings. On each sector $\Sigma_\nu:|\arg z-\nu\pi/2|<\pi/4$,
$w$ has an asymptotic expansion either (\ref{asympexw})(iv)$_b$ or else (\ref{asympexw})(iv)$_{c}$.
The Residue Theorem gives
\be{DELTA}\begin{array}{rcl}
\ds\frac1{2\pi i}\int_{\Gamma_r}w(z)\,dz&=&\ds\Delta(w)\frac{r^2}{2\pi}+o(r^2)\quad{\rm and}\cr
\ds\frac1{2\pi i}\int_{\Gamma_r}W(z)\,dz&=&-n(r,w)\quad({\rm note~that~} \Res_p W=-1),\end{array}\ee
where $\Gamma_r$ denotes the loop that was constructed in section~\ref{GAMMA}, and $\Delta(w)$ denotes the difference between the number of strings $\p$ with residue $+1$ and $-1$.
Since $|w(z)|=O(|z|)$ holds on $\Gamma_r$, and the length of the part of $\Gamma_r$ that is contained in $|\arg z-(2\nu+1)\pi/4|<\epsilon$ is at most
$4\pi\epsilon r$, the contribution of the sector $\Sigma_\nu$ to the first integral (\ref{DELTA}) is $o(r^2)$ if $w$ has the asymptotic expansion
$w\sim \pm\gamma/z+\cdots$, and is $\ds (-1)^{\nu-1}{r^2}/\pi+o(r^2)$ if $w\sim -2z+\cdots$.
In other words, the contribution of the sector $\Sigma_\nu$ to $\Delta(w)$ is $0$ if $w\sim \pm\gamma/z+\cdots$, and is $-2(-1)^{\nu}$ if $w\sim -2z+\cdots$ holds on $\Sigma_\nu$.
We apply the {\sc B\"acklund} transformation (\ref{BackIV}) to obtain
$$w_1=\frac{w'-2\gamma-2zw-w^2}{2w}$$
and conclude  from (\ref{wundW}) and (\ref{DELTA}) that
\be{CONTRIBUTION}n(r,w_1)-n(r,w)=(\Delta(w)+\Delta(w_1))\frac{r^2}{2\pi}+o(r^2).\ee
The main idea of proof now is to determine some appropriate {\sc B\"acklund} transformation such that $\Delta(w)+\Delta(w_1)$ is negative.
To this end we introduce the {\it signature} $\schema{a_0}{a_1}{a_2}{a_3}$ of $w$ to indicate that
$$W=2a_\nu z+\frac{\lambda_\nu}z+\GO{-3}$$
holds on the {\sc Stokes} sector $\Sigma_\nu$, uniformly on every closed sub-sector, with  $a_\nu\in\{-\gamma,\gamma,\alpha\}$ and $\lambda_\nu$ according to (\ref{asympexW}).
Before going into details we shall study the entire function $\ds F(z)=e^{-\int W(z)\,dz}$ with corresponding signature
$\schema{a_0}{a_1}{a_2}{a_3}$ of $w$.
Corollary~\ref{CORPL} in Appendix A, applied to
$$h(\zeta)=F(e^{\frac14(2\nu+1)\pi i}\sqrt\zeta)\quad(\sqrt\zeta>0 {\rm~if~}\zeta>0)$$
with  $a=(-1)^{\nu-1} i a_{\nu+1}$ and $b=(-1)^{\nu-1} i a_{\nu}$ then yields $\Re a=\Re b$ and $\Im a\le\Im b$, hence
\be{BEDINGUNG}\Im a_{\nu+1}=\Im a_\nu\quad{\rm and}\quad(-1)^{\nu-1}\Re a_{\nu+1}\le(-1)^{\nu-1}\Re a_\nu;\ee
in other words, {\it the differences $a_1-a_0$, $a_1-a_2$, $a_3-a_2$ and $a_3-a_0$ are real and non-negative.}
Moreover, $a_{\nu+1}=a_\nu$ implies $\lambda_{\nu+1}=\lambda_\nu$, and again from the {\sc Phragm\'en}-{\sc Lindel\"of} Principle it follows that
the half-plane ${\sf H}_\nu:|\arg z-(2\nu+1)\pi/4|<\pi/2$ is `pole-free' and thus $w$ has an asymptotic expansion on ${\sf H}_\nu$ by Proposition~\ref{ADJACENT}.

We are now looking for an appropriate starting point for the final argument.

\begin{itemize}\item[{\bf a.}] The change of variables $v(z)=\overline{w(\bar z)}$ transforms (IV) into equation (IV) for $v$ with parameters $\bar \alpha$ and $\bar\beta$
and first integral $V(z)=\overline{W(\bar z)}$, while the signature $\schema{a_0}{a_1}{a_2}{a_3}$ changes into $\schema{\bar a_0}{\bar a_3}{\bar a_2}{\bar a_1}$.
We may thus {\it formally interchange the positions of $a_1$ and $a_3$, while the positions of $a_0$ and $a_2$ remain fixed.}

\item[{\bf b.}] The change of variables $v(z)=-iw(iz)$  transforms (IV) into equation (IV) for $v$ with parameters $-\alpha$ and $\beta$, while the asymptotic expansions
$W\sim 2\alpha z$ and $W\sim \pm 2\gamma z$ are changed into $V\sim -2\alpha z$ and $V\sim \mp 2\gamma z$, respectively.
Thus the signature $\schema{a_0}{a_1}{a_2}{a_3}$ is changed into $\schema{-a_1}{-a_2}{-a_3}{-a_0}$, which means that we {\it formally may turn the signature by an
angle of $~90^\circ$ clockwise without changing our hypotheses.}\end{itemize}

\noindent Applying these transformations (which, in an abstract sense, generate the dihedral group $D_4$ acting on signatures) several times, if necessary,
it is easily seen that there are five {\it abstract} combinatorial configurations to be discussed:
\begin{center}{\bf 1.}~$\schema a a b b$, {\bf 2.}~$\schema a b a b$, {\bf 3.}~$\schema a a a b$, {\bf 4.}~$\schema a b a c$, and {\bf 5.}~$\schema a  a c b$,\end{center}
with $a,b,c\in\{\alpha,-\gamma,\gamma\}$; although $a$ and $b$, say, may be numerically equal, they re\-present different symbols.  For example, applying first {\bf b.}\ followed by {\bf a.}\ yields
$\schema b a c a \mapsto \schema{-a}{-c}{-a}{-b}\mapsto \schema{-\bar a}{-\bar b}{-\bar a}{-\bar c};$ the latter formally has the shape of the above signature {\bf 4.}
By applying Proposition~\ref{BACKCHANGE}, if necessary, we may identify $a$ with $\alpha$, and since we are free to choose the branch of $\gamma=\sqrt{-\beta/2}$,
in other words to replace $\gamma$ by $-\gamma$ in (\ref{BackIV}),
we may assume $b=-\gamma$ and $c=\gamma$, this leading to five {\it concrete} configurations as follows:

\begin{center}{\bf 1.}~ $\schema{\alpha}\alpha{-\gamma}{-\gamma}$,
{\bf 2.}~ $\schema\alpha{-\gamma}\alpha{-\gamma}$,
{\bf 3.}~ $\schema\alpha\alpha\alpha{-\gamma}$, {\bf 4.}~ $\schema\alpha{-\gamma}\alpha{\gamma}$, and {\bf 5.}~ $\schema\alpha\alpha{\gamma}{-\gamma}.$\end{center}

In the first case, $w$ has the asymptotic expansions (\ref{asympexw})(iv)$_b$ on the right upper half-plane $\Im z>-\Re z$, and  (\ref{asympexw})(iv)$^-_c$ on  $\Im z<-\Re z$.
Thus (\ref{BEDINGUNG}) applies with $\nu=1$: $a_1=\alpha$, $a_2=-\gamma,$  and $\nu=3$: $a_3=-\gamma$, $a_4=a_0=\alpha,$ this yielding $\alpha=-\gamma$. From
Corollary~\ref{CORPL} in Appendix A it then follows that $w$ has only finitely many poles, hence is a rational function in contrast to our assumption.

We will now show that none of the remaining cases {\bf 2.}\ to {\bf 5.}\ can occur for non-{\sc Weber-Hermite} sub-normal solutions.
Non-{\sc Weber-Hermite} solutions admit unrestricted application of the {\sc B\"acklund} transformation (\ref{BackIV}), this leading to sequences $(\alpha_k)$ and $(\gamma_k)$ of parameters:
$\alpha,\alpha_1,\alpha+1,\alpha_1+1,\ldots$ and $\gamma,\gamma_1,\gamma+1,\gamma_1+1,\ldots$
We have to distinguish two cases as follows:

{\bf Case a.} $\gamma_k\ne 0$ as long as $w_{k}\not\equiv 0$. Applying (\ref{BackIV}) twice we obtain in case {\bf 2.}
$$ \schema\alpha{-\gamma}\alpha{-\gamma} \mapsto\schema{\gamma_{1}}{-\gamma_{1}}{\gamma_{1}}{-\gamma_{1}}
 \mapsto\schema{\alpha_2}{-\gamma_2}{\alpha_2}{-\gamma_2},~\Delta(w)+\Delta(w_1)=-4+0=-4.$$
From (\ref{CONTRIBUTION}) (it is obvious that we may construct $\Gamma_r$ in such a way that it simultaneously works for $w$ and $w_1$) it follows that
$$n(r,w_{1})-n(r,w)=-4\frac{r^2}{2\pi}+o(r^2).$$
Repeating this process we obtain by $2$-periodicity of the sequence  $(\Delta(w_{k}))$ and the corresponding sequence of signatures
$$n(r,w_{k})=n(r,w)-4k\frac{r^2}{2\pi}+o(r^2);$$
in other words, at every step four strings of poles get lost. This, however, cannot be true for every $k\in\N$, that is, starting with some sub-normal solution $w$ with signature {\bf 2.}, there exists some $k$ such that
$$w_{k+1}=\frac{w_k'-2\gamma_k-2zw_k-w_k^2}{2w_k}\equiv 0,$$
and $w_k$ satisfies the {\sc Weber-Hermite} equation
$w_k'=2\gamma_k+2zw_{k}+w_k^2$ and has signature $\schema{\alpha_k}{-\gamma_k}{\alpha_k}{-\gamma_k}$.
The argument is similar in all other cases:
$${\bf 3.} \schema\alpha\alpha\alpha{-\gamma}\mapsto\schema{\gamma_1}{\gamma_1}{\gamma_1}{-\gamma_1}
\mapsto\schema{\alpha_2}{\alpha_2}{\alpha_2}{-\gamma_2}, ~\Delta(w)+\Delta(w_1)=-2+0=-2,$$
$${\bf 4.} \schema\alpha{-\gamma}\alpha{\gamma}\mapsto\schema{\gamma_1}{-\gamma_1}{\gamma_1}{\alpha_1}
\mapsto\schema{\alpha_2}{-\gamma_2}{\alpha_2}{\gamma_2},~\Delta(w)+\Delta(w_1)=-4+2=-2,$$
$${\bf 5.} \schema\alpha\alpha{\gamma}{-\gamma}\mapsto\schema{\gamma_1}{\gamma_1}{\alpha_1}{-\gamma_1}
\mapsto\schema{\alpha_2}{\alpha_2}{\gamma_2}{-\gamma_2}, ~\Delta(w)+\Delta(w_1)=\phantom{-}0-2=-2.$$

Thus $\ds n(r,w_k)=n(r,w)-kr^2/\pi+o(r^2)$ holds in all cases, this leading to the same conclusion as in case {\bf 2.}
{\it A posteriori} it turns out that the cases {\bf 4.}\ and {\bf 5.}\ will never occur ($\gamma\ne 0$).

{\bf Case b.} $\gamma=0$. Here we just have to consider the cases {\bf 2.} and {\bf 3.},
where the first step leads to $w_1$ with signature $\schema{\gamma_1}{-\gamma_1}{\gamma_1}{-\gamma_1}$ and $\schema{\gamma_1}{\gamma_1}{\gamma_1}{-\gamma_1}$,
respectively. Then $\gamma_1=0$ implies $w_1\equiv 0$ and we are done, while otherwise we may proceed like in {\bf Case a.}
This completes the proof of Theorem~\ref{SUBNORMALIV}. \Ende

\subsection{\it Rational and sub-normal solutions}
Rational solutions to equation (IV)$_{\alpha,\beta}$ are uniquely determined. They occur in three shapes
$w=-2z+\cdots$, $w=\gamma/z+\cdots$ ($\gamma^2=-\beta/2$), and $w=-\frac23 z+\cdots$
for certain well-known parameters, see \cite{GLS}, \S 26.
The sub-normal (={\sc Weber-Hermite}) solutions to (IV)$_{\alpha,\beta}$ form a one-parameter family and occur if and only if either $\beta=-2(1+2n-\alpha)^2$ or else $\beta=-2n^2,$ while $\alpha$ is arbitrary. For $\alpha\notin\Z$ there exist four sub-normal solutions that are truncated along two adjacent rays, while for $\alpha\in\Z$ their role is taken by rational solutions.
Any other solution satisfies $r^2=o(n(r,w))$.

\subsection{\it Sub-normal second transcendents}The analogs to Theorem~\ref{SUBNORMALIVnot} and \ref{SUBNORMALIV}
have been proved in \cite{NStSubnormal}. We will now show how the proof
of Theorem~\ref{SUBNORMALIVnot} and \ref{SUBNORMALIV} may be adapted to this case to obtain a proof that is quite different from the original proof in \cite{NStSubnormal}.\medskip

{\it Sketch of proof.} Let $w$ be any sub-normal solution ($n(r_k,w)=O(r_k^2)$) to equation (II). Then
$w$ has cluster set either $\CL(w)=\{-1/4\}$ or else $\CL(w)=\{0\}$.
The poles of $w$ are arranged in finitely many strings; each string has counting function $n(r,{\p})=\frac{\sqrt{2}}{3\pi}r^{3/2}+o(r^{3/2})$
and $n(r,{\p})=\frac{1}{3\pi}r^{3/2}+o(r^{3/2})$, and is asymptotic to some {\sc Stokes} ray $\arg z=2\nu\pi/3$ and $\arg z=(2\nu+1)\pi/3$, respectively.
On the sectors between the {\sc Stokes} rays, $w$ has asymptotic expansions (\ref{asympexw})(ii)$_b$ and (ii)$_a$, respectively.
The cluster set, the asymptotics, and the growth of the {\sc Nevanlinna} characteristic are invariant under the {\sc B\"acklund} transformations
\be{BackII}\begin{array}{lr}
B_+:~w\mapsto\tilde w=-w-\ds\frac{\alpha+1/2}{w'+w^2+z/2}&(\alpha\ne-1/2,~\alpha\mapsto\alpha+1=\tilde\alpha)\cr
B_-:~\tilde w\mapsto w=-\tilde w+\ds\frac{\tilde\alpha-1/2}{\tilde w'-\tilde w^2-z/2}&(\tilde \alpha\ne 1/2,~\tilde \alpha\mapsto\tilde\alpha-1=\alpha)
\end{array}\ee
Set $\tilde w=B_+[w]$ and observe that $w+\tilde w=\ds\frac{\alpha+1/2}{\tilde w'-\tilde w^2-z/2}=-\frac{\alpha+1/2}{w'+w^2+z/2}$, hence
$$w'+\tilde w'=-w^2+\tilde w^2=-W'+\tilde W'\quad{\rm and}\quad \tilde W-W=\tilde w+w+const.$$
This implies $\tilde W-W=\GO{1/2}$ $(z\notin\P_\delta(w)\cup\P_\delta(\tilde w))$ and $\CL(\tilde w)=\CL(w)$.
Sub-normal solutions with $\CL(w)=\{0\}$ do not exist; the proof is the same as proof of Theorem~\ref{SUBNORMALIVnot}.
For $\CL(w)=\{-1/4\}$, the number $\Delta(w)$, which denotes the difference between the number of strings with residue $1$ and $-1$, respectively, is non-zero and invariant under (\ref{BackII}); actually
$\Delta(w)=\pm 1$ or else $\Delta(w)=\pm 3$, since the contribution to $\Delta(w)$ of each sector
$2\nu\pi/3<\arg z<(2\nu+2)\pi/3$ is $1$ or $-1$.
Replacing $w$ by $-w$, if necessary, one may assume $\Delta(w)>0$, and with $\tilde w=B_+[w]$  it follows that
$$\begin{array}{rcl}
n(r,w)-n(r,\tilde w)&=&\ds\frac1{2\pi i}\int_{\Gamma_r}(\tilde W(z)-W(z))\,dz\cr
&=&\ds\frac1{2\pi i}\int_{\Gamma_r}(\tilde w(z)+w(z))\,dz=2\Delta(w)\frac{\sqrt{2}}{3\pi}r^{3/2}+o(r^{3/2}).
\end{array}$$

\begin{rem}\rm The {\sc Airy} equation $w'=z/2+w^2$
has three uniquely determined solutions with $\Delta(w)=-1$ and a single string of poles, see \cite{NStRiccati}; $k$-fold
application of the above {\sc B\"acklund} transformation  yields {\sc Airy} solutions with $2k+1$ strings in the same direction; $k+1$ strings have $\Res_\p w=-1$ and $k$ strings have $\Res_\p w=+1$, thus
$\ds n(r,w_k)\sim(2k+1){\sqrt{2}}r^{3/2}/{3\pi}$. In the generic case there are three `active' directions, $\Delta(w)=-3$ and $\ds n(r,w_k)\sim(2k+1){\sqrt{2}}r^{3/2}/{\pi}$.\end{rem}

\section{\bf The Distribution of Zeros and Poles}\label{PolesZeros}

\subsection{\it Equivalence classes and the order of sub-normal solutions} Fourth transcendents  $w_1$ and $w_2$ are called {\it equivalent}, if $w_1$ and $w_2$ are linked by some trivial {\sc B\"acklund} transformation
(any combination of rotations $w_2(z)=\eta w_1(\bar\eta z)$ with $\eta^4=1$ and reflections $w_2(z)=\overline{w_1(\bar z)}$).
Sub-normal transcendents $w$ are mapped by repeated application of {\sc B\"acklund} transformations onto some solution of one of the {\sc Weber-Hermite} equations
(\ref{WEBERHERMITE}): $w_0'=2\gamma_0\pm(2zw_0+w_0^2)$, $\gamma_0=-1\mp\alpha_0$. The smallest number of non-trivial transformations needed is called the {\it order} of $w$.
Sub-normal solutions have strings of poles either in each {\sc Stokes} direction (generic case) or else only in two consecutive directions (exceptional case); exceptional solutions are truncated along two {\sc Stokes} rays; they are uniquely determined and exist if and only if the final equation (\ref{WEBERHERMITE}) has no rational solution.
Henceforth we will restrict to the generic case, and leave the exceptional case to the interested reader.
From (\ref{ALPHAGAMMA}) and the proof of Theorem~\ref{SUBNORMALIV} we obtain:

\begin{thm}\label{EVENODD}The equivalence classes of generic sub-normal solutions of even and odd order $2k$ and $2k-1$ are represented by sub-normal solutions with signature
\be{GENSIGNeven}\schema\alpha{-\gamma}\alpha{-\gamma}~ (\alpha+\gamma=-2k-1)\ee
and
\be{GENSIGNodd}\schema{\gamma}{-\gamma}{\gamma}{-\gamma}~(\gamma=-k),\ee
respectively; $\alpha$ is arbitrary and $\beta=-2\gamma^2$. In each {\sc Stokes} direction these solutions have $k+1$ strings of poles with residue $-1$ and $k$ strings of poles with residue $1$ if the order is even, and $k$ strings of poles each with residue $\pm 1$ if the order is odd.\end{thm}

\begin{rem}\rm Any solution of order $2k$ and $2k-1$ with signatures (\ref{GENSIGNeven}) and (\ref{GENSIGNodd}), respectively, is embedded in a chain
\be{CHAIN}\begin{array}{c}\cdots\stackrel{(\ref{BackIV})}{\rightarrow}w_{2k}\stackrel{(\ref{BackIV})}{\rightarrow}w_{2k-1}\stackrel{(\ref{BackIV})}{\rightarrow}\cdots \stackrel{(\ref{BackIV})}{\rightarrow}w_{1}\stackrel{(\ref{BackIV})}{\rightarrow}w_{0}\cr
 w_0'=2\gamma_0+2zw_0+w_0^2\quad(\alpha_0+\gamma_0=-1,~\Res_p w_0=1).\end{array}\ee
We note that the case $\gamma_0=0$ does not occur, since otherwise
$$\ds w_1=-\frac{w_0'+2zw_0+w_0^2}{2w_0}=-2z-w_0$$
has order $1$ and solves $w_1'=-2+2zw_1+w_1^2$ ($\gamma_1=-1$, $\alpha_1=0$), in contrast to the definition of the order.\end{rem}

\subsection{\it The distribution of residues}Let $w$ be any sub-normal fourth transcendent. With each string of poles $\p=(p_k)$ we will associate a polygon $\pi(\p)$ with vertices $p_k$. These polygons divide $|z|>R$ (sufficiently large)
into domains $D_\nu$, arranged in cyclic order; $\pi_\nu=\pi(\p_\nu)$ separates $D_\nu$ from $D_{\nu-1}$. Also $D_\nu\cap\{z:|z|=r\}$ ($r>R$) is an arc of angular measure $\Theta_\nu(r)$, where
$\Theta_\nu(r)$ either tends to $\frac 34\pi$, $\frac\pi 4$ or else $0$ as $r\to\infty$; in any case $r\Theta_\nu(r)$ tends to infinity. The polygon $\pi_\nu$
is accompanied by polygons $\hat\pi_{\nu-1}$ and $\tilde\pi_\nu$ in  $D_{\nu-1}$ and $D_\nu$, respectively, which start
at $|z|=R$ such that  $\arg z\to k\frac\pi 4$ with $k=k_\nu$ as $|z|\to\infty$ on $\hat\pi_{\nu-1}$ and $\tilde\pi_\nu$, and the angular measure and the length of the shorter arc on $|z|=r$ joining these polygons to $\pi_\nu$ tends to zero and infinity, respectively, as $r\to\infty$.
Re-scaling along any sequence $(h_n)$ on $D_{\nu}$ with $|h_n|\dist(h_n,\partial D_{\nu})\to\infty$ yields a constant limit function $\w\equiv\tau_{\nu}\in\{0,-2\}$; in particular, $w=\tau_{\nu} z+\ko{}$ holds as $z\to\infty$ on $\tilde\pi_{\nu}$, and similarly we obtain $w=\tau_{\nu-1} z+\ko{}$
on $\hat\pi_{\nu-1}$. To determine $\epsilon_\nu=\Res_{\p_\nu}w\in\{-1,1\}$ we assume for simplicity $\arg z\sim \pi/4$ on $\pi_\nu$, and compute
$$\frac1{2\pi i}\int_{\kappa_\nu(r)} w(z)\,dz=\epsilon_\nu\frac{r^2}{2\pi}+o(r^2)$$
along the positively oriented simple closed curve $\kappa_\nu(r)$ which
consists of sub-arcs $\hat\kappa_{\nu-1}(r)$ and $\tilde\kappa_\nu(r)$ of $\hat\pi_{\nu-1}$ and $\tilde\pi_\nu$ joining $|z|=r_0$ to $|z|=r$, respectively,
and sub-arcs $\sigma_{r_0}$ of $|z|=r_0>R$ and $\sigma_r$ of $|z|=r>r_0$. The latter has length $o(r)$, hence
$\frac1{2\pi i}\int_{\sigma_{r}}w(z)\,dz =o(r^2)$ holds. From $w(z)=\tau_{\nu-1}z+o(|z|)$ and $w(z)=\tau_{\nu}z+o(|z|)$  on $\hat\pi_{\nu-1}$  and $\tilde\pi_\nu$, respectively, and $\arg z\to\pi/4$ as $|z|\to\infty$ we obtain
$$\frac1{2\pi i}\int_{\kappa_{\nu}(r)}w(z)\,dz=\frac{(\tau_{\nu-1}-\tau_\nu)r^2}{4\pi}+o(r^2),$$
which implies $\ds\epsilon_\nu=(\tau_{\nu-1}-\tau_\nu)/{2}$. The relation between the asymptotics and the residues in different {\sc Stokes} directions is displayed in Table 1. \medskip

\begin{center}\begin{tabular}{|c|c|c|c|}\hline
{\sc Stokes} ray&$\tau_{\nu-1}$&$\tau_\nu$&$\epsilon_\nu$\cr\hline
$\arg z=\frac\pi 4, \frac54\pi$&$\begin{array}{c}\phantom{-}0\cr -2\end{array}$&$\begin{array}{c}-2\cr\phantom{-}0\end{array}$&$\begin{array}{c}\phantom{-}1\cr- 1\end{array}$\cr\hline
$\arg z=\frac34\pi, \frac74\pi$&$\begin{array}{c}\phantom{-}0\cr -2\end{array}$&$\begin{array}{c}-2\cr \phantom{-}0\end{array}$&$\begin{array}{c}-1\cr \phantom{-}1\end{array}$\cr\hline\end{tabular}

\medskip \small{\sc Table 1.} Asymptotics and distribution of residues.\end{center}

\begin{ex}\rm  Any sub-normal solution $w$ with signature $\schema\alpha{-\gamma}\alpha{-\gamma}$ (and parameter $\gamma=-(2k+1+\alpha)$)
has $2k+1$ strings of poles asymptotic to  $\arg z=\pi/4$; in counter-clockwise order the residues are
 $(-1,+1,-1,\ldots,+1,-1).$ If $w$ has signature $\schema k{-k}k{-k}$,
the number of strings is $2k$, and the residue
vector in the same direction is given by $(+1,-1,\ldots,+1,-1)$, see also Table~2 below.\end{ex}

\subsection{\it The distribution of zeros} The re-scaling method does not immediately apply to detect the zeros of $w$, since re-scaling of any fourth transcendent along any sequence of zeros yields the limit function $\w\equiv 0$. The reason for this is that $w'=\pm 2\gamma$ at zeros
is `small', hence the initial values for the limit function are $\w(0)=\w'(0)=0$. Nevertheless it is possible to determine the distribution of zeros in any case, and, in particular, for sub-normal solutions. The zeros of $w$ with $w'=\pm 2\gamma$ are poles of $w_\pm=\ds\frac{w'\pm 2\gamma}{2w}+z+\frac{w}{2}$ with residue $1$. We thus may conclude that the zeros
are distributed in the same manner as are the poles. In particular, they form strings if $w$ is sub-normal.\medskip

\begin{center}{\small$\begin{array}{|cccccc|}\hline
&&&&\star&\circ\cr
&&&\star&\circ&\bullet\cr
&&\star&\circ&\bullet&\doppel\cr
&\star&\circ&\bullet&\doppel&\circ\cr
\star&\circ&\bullet&\doppel&\circ&\cr
\circ&\bullet&\doppel&\circ&&\cr
\bullet&\doppel&\circ&&&\cr
\doppel&\circ&&&&\cr\hline
\end{array}$\quad
$\begin{array}{|cccccc|}\hline
&&&&\doppel&\circ\cr
&&&\doppel&\circ&\bullet\cr
&&\doppel&\circ&\bullet&\doppel\cr
&\doppel&\circ&\bullet&\doppel&\circ\cr
\doppel&\circ&\bullet&\doppel&\circ&\bullet\cr
\circ&\bullet&\doppel&\circ&\bullet&\cr
\bullet&\doppel&\circ&\bullet&&\cr
\doppel&\circ&\bullet&&&\cr\hline
\end{array}$}

\medskip {\small {\sc Table 2.} Distribution of poles $\bullet$~$\circ$ ($\Res_\bullet w=1$, $\Res_\circ w=-1$) and zeros $\star$$\ast$ along $\arg z=\pi/4$;
$w$ has signature (\ref{GENSIGNeven}) and (\ref{GENSIGNodd}), and order $2$ and $3$, respectively.
In case of $\gamma=0$ the double string of zeros $\doppel$ on the left hand side collapses to a string of double zeros, while the single string $\star$ disappears.
The deficiency of zero then is $\delta(0,w)=1/3$.}
\end{center}

\subsection{\it {Painlev\'e} transcendents and first order differential equations}\label{FIRSTORDERDGL}
It is obvious that every {\sc Airy}- and {\sc Weber-Hermite} solution also satisfies some first order algebraic differential equation
\be{FIRSTORDER}P(z,w,w')={w'}^n+\sum_{\nu=1}^{n-1}P_\nu(z,w){w'}^\nu=0;\ee
$P_\nu$ is a polynomial in $w$ (of degree $\le 2n-2\nu$) over the field of rational functions in $z$. The converse was proved in \cite{GLS}, Theorem~21.1
and Theorem~25.4 for second and fourth transcendents, respectively.
Based on the re-scaling method we will give a quite different proof of

\begin{thm}\label{1Order} The second and fourth {\sc Painlev\'e} transcendents also satisfying some
first order algebraic differential equation
{\rm (\ref{FIRSTORDER})} coincide with the {\sc Airy}- and {\sc Weber-Hermite} solutions, respectively, while first {\sc Painlev\'e} transcendents never solve first order equations.\end{thm}

\beweis We consider solutions to (IV) which also solve (\ref{FIRSTORDER}) and assume that
$$P(z,x,y)=\prod_{\nu=1}^n(y-G_\nu(z,x))$$
is {\it irreducible}. From
$P(h+h^{-1}\z,hw_h(\z),h^2w'_h(\z))=0$ it then follows that the limit functions of the re-scaling process satisfy some {\sc Briot-Bouquet} differential equation
$$Q(\w,\w')=0$$
with $Q(\x,\y)=\lim_{h\to\infty}h^{-2n}P(h,h\x,h^2\y)=\lim_{h\to\infty}h^{-2n}\prod_{\nu=1}^n(h^{2}\y-G_\nu(h,h\x))$.
On the other hand it is known that at any pole $p\ne 0$ with residue $\epsilon$,
$$\epsilon w'+2zw+w^2=O(|p|^2)\quad(p\to\infty)$$
holds on $|z-p|<\delta|p|^{-1}$. This implies
$G_\nu(z,x)=\epsilon_\nu x^2+2\epsilon_\nu zx+O(|z|^2)$ as $x\to\infty$,
uniformly with respect to $z$, $|z-p|<\delta|p|^{-1}$. Re-scaling about any sequence of poles yields
$$Q(\w,\w')=\prod_{\nu=1}^n(\w'-\epsilon_\nu(\w^2+2\w+\a_\nu))=0\quad(\epsilon_\nu=\pm 1,~\a_\nu\in\C).$$
Thus $\w$ satisfies $\w'=\epsilon(\a+2\w+\w^2)$ for some $\epsilon\in\{1,-1\}$ and $\a\in\C$, which is compatible with (\ref{PainleveElliptic})(iv) if and only if $\a=\c=0$.
In particular, it follows that $w$ has cluster set $\CL(w)=\{0\}$, hence the set $\P$ of poles has {\it string structure}, that is, $\P$ consists of finitely or infinitely many strings of poles.
In the first case $w$ is sub-normal, hence a {\sc Weber-Hermite} solution. To rule out the second possibility we have to discuss two subcases, which occur at every first order differential equation (see {\sc Eremenko}~\cite{Eremenko, Eremenko1}): equation
(\ref{FIRSTORDER}) has {\it genus} either $g=0$ or else
$g=1.$ In the first case we have $w=R(z,y)$, where $R$ is rational and $y$ satisfies some {\sc Riccati} equation with rational coefficients. Solutions to {\sc Riccati} equations have only finitely many strings of poles, and this also holds for $w$ itself. In the second case,
$w$ is a rational function of $u$ and $u'$ over the field of algebraic functions, where $u$ satisfies some differential equation
$$u'^2=4a(z)(w-e_1)(w-e_2)(w-e_3)$$
with $e_\mu\ne e_\nu$ for $\mu\ne\nu$, $e_1+e_2+e_3=0$, and $a$ an algebraic function (thus $u(z)=\wp(\int\sqrt{a(z)}\,dz)$ locally).
This also leads to a contradiction, since this time the set $\P$ of poles locally has {\it lattice structure} in contrast to its string structure.

The proof for second transcendents runs along the same way.
It remains to prove that first transcendents do not solve any first order differential equation (Theorem~13.1 in \cite{GLS}).
This time we obtain the {\sc Briot-Bouquet} differential equation
$$Q(\w,\w')=\prod_{\nu=1}^n(\w'-2\epsilon_\nu\w^{3/2}-\a_\nu)=0\quad(\epsilon_\nu=\pm 1,~\a_\nu\in\C),$$
in contrast to $\w'^2=4\w^3+2\w-2\c$. \Ende

\begin{ex}\rm To obtain first order equations for sub-normal solutions start with $w_0'=2\gamma_0+2zw_0+w_0^2$ ($\alpha_0=-1-\gamma_0$), say,
and compute successively (and in a purely algebraic manner)
$$\ds w_{\nu+1}=-\frac{w'_\nu+2\gamma_\nu+2zw_\nu+w_\nu^2}{2w_\nu}=\frac{P_\nu(z,w_0)}{Q_\nu(z,w_0)}\quad{\rm and}\quad w_{\nu+1}'=\frac{\tilde P_\nu(z,w_0)}{\tilde Q_\nu(z,w_0)}.$$
The {\it resultant} with respect to $w_0$ of the polynomials $w_{\nu+1}Q_\nu(z,w_0)-P_\nu(z,w_0)$ and $w'_{\nu+1}\tilde Q_\nu(z,w_0)-\tilde P_\nu(z,w_0)$
is a polynomial $R_{\nu+1}(z,w_{\nu+1},w_{\nu+1}')$. Then $w=w_k$ is sub-normal of order $k$ and satisfies $R_k(z,w,w')=0$; the parameter $\alpha=\alpha_k$
may be prescribed by adjusting $\gamma_0$, while $\gamma=\gamma_k$ then is fixed. In the first step we obtain
$$w'^2+4w'-w^4-4zw^3-4(z^2-\alpha)w^2+4=0\quad(\gamma=-1)$$
and the {\it binomial differential equation}
$y'^2=(y^2-4\alpha)(y-2z)^2$ for $y=w+2z.$ The re-scaling process yields $\y'^2=\y^2(\y-2)^2$ and $\w'^2=\w^2(\w+2)^2$.\end{ex}

\section{\bf Deficient Values and Functions}\label{VALUEDISTRIBUTION}
\subsection{\it The deficiency of zero of fourth transcendents}It is well known and easy to prove that
$$\ds m\Big(r,\frac1{w-c}\Big)=O(\log r)\quad(c\ne 0)$$
holds for every {\sc Painlev\'e} transcendent.
This is also true for $c=0$ except in case of equation (II) with $\alpha=0$ and equation (IV) with $\beta=0$. In \cite{NStSubnormal} it was shown that
the {\sc Nevanlinna} {\it deficiency}
$$\delta(0,w)=\liminf_{r\to\infty}\frac{m(r,1/w)}{T(r,w)}$$
vanishes for every transcendental solutions to $w''=zw+2w^3$ $(\alpha=0)$. In case of equation (IV) it is known that
either $0$ is a {\sc Picard} value and $w$ solves some {\sc Weber-Hermite} equation $w'=\pm(2zw+w^2)$, or else $\delta(0,w)\le 1/2$ holds (see \cite{GLS}).
The whole truth, however, is more refined and unlikely.

\begin{thm}\label{Nulldefekt}{\rm (also \cite{NStHamilton})}
Any transcendental solution to equation (IV) with $\beta=0$
has {\sc Valiron} deficiency
$$\ds\Delta(0,w)=\limsup_{r\to\infty}\frac{m(r,1/w)}{T(r,w)}=0,$$
except in the following case: $w$ is sub-normal of order $2k$ with
parameters $\beta=0$ and $\alpha=\pm(2k+1)$, and
$$\ds\Delta(0,w)=\delta(0,w)=\frac1{2k+1}.$$
\end{thm}\medskip

\beweis We first consider the case when $w$ is not sub-normal. Then the
simple closed curve $\Gamma_r$ that was constructed in section~\ref{GAMMA} will be used to compute
$$n_+(r,w)-n_-(r,w)=\frac1{2\pi i}\int_{\Gamma_r}w(z)\,dz=O(r^2);$$
here $n_\pm(r,w)$ denotes the number of poles on $|z|<r$ with residue $\pm 1$. Then  $n_+(r,\tilde w)-n_-(r,\tilde w)=O(r^2)$ also holds for
\be{w1w}\ds \tilde w=\frac{w'-2zw-w^2}{2w},\ee
and
$n(r,1/w)=2n_+(r,\tilde w)$ holds--note that $w$ has only double zeros, $\tilde w$ is regular at poles of $w$ with residue $-1$ and has poles with residue $-1$ at poles of $w$ with residue $1$.
This implies $n(r,1/w)=2n_-(r,\tilde w)+O(r^2)=2n_+(r,w)+O(r^2)=n(r,w)+O(r^2)$ and
$N(r,1/w)=N(r,w)+O(r^2)=T(r,w)+O(r^2).$
Thus $\lim_{r\to\infty}{T(r,w)}/{r^2}=\infty$ implies $\Delta(0,w)=\delta(0,w)=0$.
It remains to consider sub-normal solutions to equation (IV) with  signature $\schema\alpha 0\alpha 0$ (generic case);
the exceptional case $\schema\alpha\alpha\alpha 0$is dealt with in the same manner. These solutions have even order $2k$ and parameter $\alpha=-2k-1$.
Solutions of order $0$ have {\sc Picard} value zero. If, however, the order is $2k\ge 2$ we
again consider $\tilde w$ defined by (\ref{w1w}), and remind the reader that again
the (double) zeros of $w$ are poles of $\tilde w$ with residue $+1$, poles of $w$ with residue $+1$ are poles of $\tilde w$ with residue $-1$, and
$\tilde w$ is regular at all other points of the plane, including the poles of $w$ with residue $-1$.
Then $w$ and $\tilde w$ have counting functions
$n(r,w)\sim 4(2k+1)r^2/2\pi$ and $n_+(r,\tilde w)\sim n_-(r,\tilde w)\sim 4kr^2/2\pi$, respectively, where $\sim$ means up to some term $o(r^2)$.
This implies $n(r,1/w)=2n_+(r,\tilde w)\sim 8kr^2/2\pi$ and $\Delta(0,w)=\delta(0,w)=1/(2k+1)$. \Ende

\begin{rem}\label{Computem}\rm There is a second way to prove $m(r,1/w)\sim 2r^2/\pi$
for $w$ with signature $\schema\alpha 0\alpha 0$ as follows. Proceeding like in the proof of Theorem~\ref{EXUNIQUE} with $\gamma=0,$ $w_1=w$, $w_2=0$, and $y=z^{1/2}(w_1-w_2)=z^{1/2}w$ we obtain
$$y''-(4z^2-4\alpha-2+\GO{-2})y=0$$
on the {\sc Stokes} sectors $\Sigma_1$ and $\Sigma_2$, hence $-\log|w(z)|=-\Re z^2+o(|z|^2)$
on these sectors, while $-\log|w(z)|=O(\log |z|)$ holds on $\Sigma_0$ and $\Sigma_2$. This gives the assertion, details are left to the reader.\end{rem}

\subsection{\it Deficient rational functions of fourth transcendents}Suppose that $w$ denotes any meromorphic non-rational solution to some algebraic differential equation
$$\Omega(z,w,w',\ldots,w^{(n)})=0\quad(\Omega {\rm ~some~polynomial}),$$
and let $\phi$ be any rational function such that $\Omega(z,\phi,\phi',\ldots,\phi^{(n)})\not\equiv 0$. Then
$$m\Big(r,\frac{1}{w-\phi}\Big)=S(r,w),$$
and the right hand side is $O(\log r)$ if $w$ has finite order of growth. This is a special case of a theorem of {\sc Mokhonko}-{\sc Mokhonko} \cite{MokMok},
which was used before by {\sc Wittich} (see his book~\cite{HW1}, for example) in many particular cases; if $\phi$ is not rational, the term $O(T(r,\phi))+S(r,\phi)$ has to be added
.on the right hand side

To compute $\delta(\phi,w)$ for sub-normal fourth transcendents $w$ and rational  functions it suffices to consider solutions $\phi$ to the same equation.
Let $\tilde w=B[w]$ be any {\sc B\"acklund} transform of $w$. For definiteness we choose
$$B[w]=\ds\frac{w'-2\gamma-2zw-w^2}{2w}$$
and set $\tilde\phi=B[\phi]$. Then $2w\tilde w-2\phi\tilde\phi=w'-\phi'-2z(w-\phi)-(w^2-\phi^2)$, hence
$$2(w-\phi)\tilde w+2\phi(\tilde w-\tilde\phi)=w'-\phi'-(w-\phi)(2z+w+\phi)$$
holds, and dividing by $w-\phi$ yields
$\ds 2\phi\frac{\tilde w-\tilde\phi}{w-\phi}=-2\tilde w+\frac{w'-\phi'}{w-\phi}-2z-w-\phi.$
By the usual rules we thus obtain
$\ds m\Big(r,\frac{\tilde w-\tilde \phi}{w-\phi}\Big)=O(\log r)$,
and $\ds m\Big(r,\frac{ w-\phi}{\tilde w-\tilde\phi}\Big)=O(\log r)$ by symmetry, thus
\be{mtilde} m\Big(r,\frac{1}{w-\phi}\Big)= m\Big(r,\frac1{\tilde w-\tilde \phi}\Big)+O(\log r).\ee
The argument works as well for any other and also for iterated {\sc B\"acklund} transformations. Now every rational solution $\phi(z)=-2z+\cdots$ and $\phi(z)=\frac\gamma z+\cdots$ has its origin
in $\phi_0(z)=0$, that is, there exists some (iterated) B\"acklund transformation, again denoted $B$, such that $B[\phi]\equiv 0$. If $w$ is not sub-normal this yields
\be{mrwphi}m\Big(r,\frac{1}{w-\phi}\Big)=m\Big(r,\frac{1}{B[w]}\Big)+O(\log r)=o(T(r,w)).\ee
If, however, $w$ is sub-normal of order $k$ we take $B$ such that $\tilde w=B[w]$ satisfies some {\sc Weber-Hermite} equation, which then is also solved by $\tilde\phi=B[\phi]$.
Then $\tilde w-\tilde\phi$ has {\sc Picard} value zero, hence
$$m\Big(r,\frac1{w-\phi}\Big)=m\Big(r,\frac1{\tilde w-\tilde\phi}\Big)+O(\log r)=T(r,\tilde w)+o(r^2)=4\frac{r^2}{2\pi}+o(r^2)$$
and $\ds T(r,w)=4(k+1)\frac{r^2}{2\pi}+o(r^2)$
hold in the generic case, this implying $\Delta(\phi,w)=\delta(\phi,w)=1/(k+1)$. In the exceptional cases the factor $4$ has to be replaced by $2$, with the same result.
We note that (\ref{mrwphi}) remains true if $w$ is not sub-normal, the rational function $\phi$ is replaced by any sub-normal solution $\phi$, and the term $O(\log r)$ is replaced by $O(r^2)$. We thus have proved

\begin{thm}\label{Ratdefekt}Let $w$ be any fourth transcendent and $\phi$ be any sub-normal or rational solution $($but $\phi(z)\not\equiv -\frac23 z+\cdots)$ to the very same equation. Then
$$\Delta(\phi,w)=\delta(\phi,w)=\left\{\begin{array}{cl}0&(w {\rm ~not~sub\!\!-\!\!normal})\cr
\ds\frac1{k+1}&(w {\rm ~sub\!\!-\!\!normal~of~order~} k\ge 1 {\rm ~and~} \phi {\rm ~rational}).\end{array}\right.$$
\end{thm}

\begin{ex}\rm The sub-normal solutions to equation (IV) with $\alpha=\beta=-2$  have order $1$, solve
$w'^2+4w'-w^4-4zw^3-4(z^2+2)w^2+4=0$, and have deficiency $\delta(\frac 1z,w)=\frac12.$ \end{ex}

\begin{rem}\rm {\sc Shimomura}~\cite{ShimomuraB} proved $\delta(\phi,w)\le 1/2$ if $\beta\ne 0$ and $\delta(\phi,w)\le 3/4$ if $\beta=0$ for fourth transcendents $w$
and so-called `small' functions satisfying $T(r,\phi)=S(r,w)$, with the obvious exception that $w$ solves some {\sc Weber-Hermite} equation.
\end{rem}

\begin{rem}\rm Rational solutions $\phi(z)=-\frac23z+\cdots$ arise for parameters $\beta=-2(2n+\frac13-\alpha)^2$ ($n\in\Z$),
$\alpha\in\Z$, and are not related to  {\sc Weber-Hermite} solutions; the most simple case is $\phi(z)=-\frac23 z$, $\alpha=0$, $\beta=-\frac29$.
The first part of the proof also works in the present case and again yields (\ref{mtilde}).
Thus the general case may be reduced to solutions $w$ to
$$\ts 2ww''=w'^2+3w^4+8zw^3+4z^2w^2-\frac49\quad{\rm and}\quad \phi(z)=-\frac23 z,$$
($\alpha=0$, $\beta=-2/9$), but nevertheless requirers a new idea. \end{rem}

\subsection{\it Second transcendents} Theorem~\ref{Nulldefekt} and \ref{Ratdefekt} have an analog for second transcendents.
Here rational and sub-normal solutions are separated from each other since they correspond to parameters $\alpha\in\Z$ and $\alpha\in\frac12+\Z$, respectively.

\begin{thm}Let $w$ be any normal second transcendent. Then
$$m\Big(r,\frac1{w-\phi}\Big)=O(r^{3/2})=o(T(r,w))$$
holds for every meromorphic function satisfying $T(r,\phi)=O(r^{3/2}).$
\end{thm}

\beweis We note that for $\alpha=0$,
\be{mr1wII}m\Big(r,\frac1{w}\Big)=O\big(\sqrt{T(r,w)}\;\big)=O(r^{3/2})\ee
has been proved in \cite{NStSubnormal}. To prove the general result we may assume that $w$ and $\phi$ satisfy one and the same {\sc Painlev\'e} equation (II).
Then either $\alpha\in\Z$ and $\phi$ is rational or else $\alpha\in\frac 12+\Z$ and $\phi$ is sub-normal.
This time the {\sc B\"acklund} transformations (\ref{BackII}) and the special transformation
\be{SpecBaecktrans}z=-\sqrt[3]{2}\,t,~-\sqrt[3]{2}\,y(t)^2=w'(z)-w(z)^2-z/2\quad(\alpha=1/2)\ee
play an important role; the latter transforms $w''=1/2+zw+2w^3$ into $\ddot y=ty+2y^3$, hence forms the bridge between $\alpha=1/2$ and $\alpha=0$, see \cite{GLS}, p.\ 142.

Consider $\ds\tilde w=B[w]=-w+\frac{\alpha-1/2}{w'-w^2-z/2}$ with $\alpha\ne 1/2$, say. Then from
$$m\Big(r,\frac1{w'-w^2-z/2}\Big)=m(r,\tilde w+w)+O(1)=O(\log r),$$
the same estimate with $w$ and $\tilde w$ replaced by $\phi$ and $\tilde\phi=B[\phi]$, and
$$\ds\tilde w-\tilde\phi=-(w-\phi)+(\alpha-1/2)\frac{\phi'-w'+(w-\phi)(w+\phi)}{(w'-w^2-z/2)(\phi'-\phi^2-z/2)}$$
it follows that $\ds m\Big(r,\frac{\tilde w-\tilde \phi}{w-\phi}\Big)=O(\log r),$
hence  (\ref{mtilde}) again holds for every (iterated) {\sc B\"acklund} transform $\tilde w=B[w]$. Like in the proof of Theorem~\ref{Ratdefekt} this enables us to reduce the general case
$\alpha\in\Z$ to the case  $\alpha=0$, hence the assertion follows from the corresponding estimate (\ref{mr1wII}).
If, however, $\alpha\in\frac12+\Z$ we may assume $\alpha=\frac12$ and $\phi'=\phi^2+z/2$, and use (\ref{SpecBaecktrans}):
Set $v=w'-w^2-z/2$; then the assertion follows from $v=w'-\phi'-(w-\phi)(w+\phi)$, hence
$$\frac1{w-\phi}=\frac1v\Big(\frac{w'-\phi'}{w-\phi}-w-\phi\Big),$$
and $m(r,1/v)=O\big(\sqrt{T(r,v)}\;\big)=O(r^{3/2})$. \Ende

\section{\bf Appendix A: The Phragm\'en-Lindel\"of Principle}\label{AppendixB}

The {\sc Phragm\'en}-{\sc Lindel\"of} Principle is an easy consequence of the Two-Constants-Theorem and may be stated as follows
(see also {\sc Titchmarsh}~\cite{ECT}, p.~176-180):
\begin{thm}\label{PhragLind}Let $f$ be a holomorphic function of finite order
$$\ds\limsup_{z\to\infty}\frac{\log^+\log^+|f(z)|}{\log|z|}<\frac\pi{2\delta}\quad{\rm on~}
\Sigma:|\arg z|\le\delta,~|z|\ge r_0.$$
Then  $f$ is bounded on $\Sigma$ if $f$ is bounded on $\partial\Sigma\setminus\{\infty\}$,
and the following is true:
\begin{itemize}
 \item[{\bf 1.}] If $f$ tends to $c$ as
$z=re^{i\delta}\to\infty$, then $f$ tends to $c$ as $z\to\infty$,
uniformly on every sector $-\delta+\epsilon<\arg z\le\delta$
$(0<\epsilon<\delta$ arbitrary$)$.
\item[{\bf 2.}] If $f$ tends to $c_{\pm}$ as $z=re^{\pm i\delta}\to\infty$, then
$c_+=c_{-}=c$, and $f$ tends to $c$ as $z\to\infty$, uniformly on $|\arg
z|\le\delta.$
 \end{itemize}\end{thm}

\begin{cor}\label{CORPL}Let $h$ be any holomorphic function of finite order
on $|\arg z|\le\eta$, $|z|\ge r_0$,  and assume that for every $\delta$, $0<\delta\le\eta$, $h$ satisfies
$$h(z)=\left\{\begin{array}{ll}C_+e^{az}z^\lambda(1+o(1))&(z=re^{i\delta}\to\infty)\cr
C_-e^{bz}z^\mu(1+o(1))&(z=re^{-i\delta}\to\infty)\end{array}\right.\quad{\rm with}~C_+C_-\ne 0.$$
Then the following is true:
\begin{itemize}\item[{\bf 1.}] $\Re a=\Re b$ and $\Im a\le \Im b$.
\item[{\bf 2.}] $a=b$ implies $\lambda=\mu$ and
$\ds h(z)=Ce^{az}z^\lambda(1+o(1))$ as $z\to\infty,$
uniformly on $|\arg z|\le\eta$; in particular, $h$ has only finitely many zeros.
\end{itemize}\end{cor}

\beweis In order to prove {\bf 1.} we assume $\Re a\ne Re b$ and even $\Re a>\Re b$
(otherwise take $\overline{h(\bar z)}$ instead of $h$ to replace
$a,b$ by $\bar b, \bar a$), and consider the holomorphic function
 \be{htof}f(z)=e^{-az}z^{-\lambda}h(z);\ee
it has finite order on $|\arg z|\le\eta$ and satisfies $f(re^{i\delta})\to
C_+\ne 0$  and
 \be{unten}|f(re^{-i\delta})|=O(r^{\Re(\mu-\lambda)})\;e^{\Re(b-a)r\cos\delta+\Im(b-a)r\sin\delta}\ee
as $r\to\infty$. The right hand side of (\ref{unten}) tends to zero
if $\delta>0$ is chosen sufficiently small, hence $f$ is bounded on $|z|\ge r_0$, $|\arg z|\le\delta$
by the {\sc Phragm\'en}-{\sc Lindel\"of} Principle. The second part of the {\sc Phragm\'en}-{\sc Lindel\"of} Principle then gives a contradiction
with $c_{+}=C_+\ne 0$ and $c_{-}=0$. In the same way obtain a
contradiction if we assume $\Re a=\Re b$ and $\Im a>\Im b:$ $f$
tends to $C_+$ as $z=re^{i\delta}\to\infty$, while (\ref{unten}) again
implies that $f(re^{-i\delta})\to 0$ as $r\to\infty$. This proves $\Im a\le\Im b$.

In order to prove {\bf 2.} we assume $a=b$ and $\Re\lambda>\Re\mu$. Then $f$, again defined by (\ref{htof}), tends to $C_+\ne 0$ as $z=re^{i\delta}\to\infty$, while
$|f(re^{-i\delta})|=O\big(r^{\Re(\mu-\lambda)}\big)$ tends to zero as $r\to\infty.$ This proves $\Re\lambda\le\Re\mu$, and in the
same way we obtain $\Re\lambda\ge\Re\mu$, hence
$\Re\lambda=\Re\mu$. This eventually implies that $f$ tends to
$C_+$ as $z\to\infty$ on $-\delta+\epsilon<\arg z\le\delta$ by the
{\sc Phragm\'en}-{\sc Lindel\"of} Principle, hence, in particular,
$h(x)= C_+e^{ax}x^\lambda(1+o(1))$ holds as $x\to +\infty$ ($x$ real). In the same way we
obtain $h(x)=C_-e^{ax}x^\mu(1+o(1))$, hence $\lambda=\mu$, $C_+=C_-=C$,
and $h(z)=Ce^{az}z^\lambda(1+o(1))$ holds as $z\to\infty$ on $|\arg z|\le\eta$. \Ende

\section{\bf Appendix B: Asymptotic Expansions}\label{AppendixC}

\subsection{\it Asymptotic expansions of specific solutions}
The following theorem on the existence of asymptotic expansions of solutions to algebraic differential equations applies in many different situations to {\it specific} solutions.

\begin{thm}\label{AsympEx}Let $w$ be any solution to the algebraic differential equation
$$Q[w]=P(z,w)$$
satisfying $w(z)\to c_0$ as $z\to\infty$ on some sector $S$, and assume that
\begin{itemize}\item[{\bf a.}] $P(z,w)$ is a polynomial in  $w$ and  rational in $z$ satisfying
$$P(z,c_0)\to 0\quad{\rm and}\quad P_w(z,c_0)\to c\ne 0\quad(z\to\infty {\rm~ on~} S);$$
\item[{\bf b.}] $Q[w]=\sum_{M}a_M(z)w^{\ell_0}M[w]$
is a differential polynomial with rational coefficients
$a_M(z)=A_Mz^{\alpha_M}(1+o(1))$ as $z\to\infty$ $(A_M\ne 0,~\alpha_M\in\Z),$
and monomials $M[w]=w'^{\ell_1}\cdots{w^{(m)}}^{\ell_m}$ of weight
$$d_M=2\ell_1+\cdots +(m+1)\ell_m\ge \alpha_M+2.$$
\end{itemize}
Then $w$ has an asymptotic expansion
$w\sim \sum\limits_{k=0}^\infty c_kz^{-k}$ on $S.$\end{thm}

\beweis We will start with
 \be{ASYMw_n}w(z)=\sum_{\nu=0}^nc_\nu z^{-\nu}+o(|z|^{-n})=\psi_n(z)+o(|z|^{-n}),\ee
which is true for $n=0$. To proceed further we need
\be{HPSI}w^{\ell_0}M[w]=\psi_n^{\ell_0}M[\psi_n]+o(|z|^{-d_M-n+1})\quad(z\to\infty);\ee
the proof will be given below. Using (\ref{HPSI}) we obtain
$$Q[w]=Q[\psi_n]+o(|z|^{\max_M(\alpha_M-d_M)-n+1})=Q[\psi_n]+\ko{-n-1},$$
hence $w$ satisfies
\be{ALG_w}P(z,w)=Q[\psi_n]+\ko{-n-1}.\ee
It follows from the first hypothesis that for any rational function $R$ that tends to zero as $z\to\infty$, the algebraic equation $P(z,y)=R(z)$ has a unique solution
that tends to $c_0$ as $z\to\infty$. In particular, equation
\be{ALG_a}P(z,y)=Q[\psi_n](z)\ee
has a unique solution $\ds y_n(z)=\sum_{\nu=0}^\infty a_\nu^{[n]}z^{-\nu}$
about $z=\infty$, and from (\ref{ALG_w}), (\ref{ALG_a}), and
$$P(z,w(z))-P(z,y_n(z))=\int_{y_n(z)}^{w(z)}P_\zeta(z,\zeta)\,d\zeta
=(c+o(1))(w(z)-y_n(z))$$
as $z\to\infty$ (we integrate along the straight line in $|\zeta-c_0|<\delta$ from $y_n(z)$ to $w(z)$) it follows that
 $$w(z)-y_n(z)=o(|z|^{-n-1}),$$
hence (\ref{ASYMw_n}) holds with $n$ replaced by $n+1$ and $\psi_{n+1}(z)=\sum\limits_{\nu=0}^{n+1}a^{[n]}_\nu z^{-\nu}$; we note that  $a_\nu^{[k]}=c_\nu$ ($0\le\nu\le k\le n$) holds, while the new coefficient is $c_{n+1}=a_{n+1}^{[n]}.$

To prove (\ref{HPSI}) we first consider the case $\psi_n(z)\equiv c_0$, hence $M[\psi_n]=0$,
$w^{(k)}(z)=\psi_n^{(k)}(z)+o(|z|^{-n-k})=o(|z|^{-n-k}),$ and
$$w^{\ell_0}M[w]=o(|z|^{\ell_1(-n-1)+\cdots+ l_m(-n-m)})=o(|z|^{-d_M-n+1})$$
since $-\sum_{k=1}^m(n+k)\ell_k\le -d_M-(n-1)\sum_{k=1}^m\ell_k\le-d_M-n+1.$ Thus (\ref{HPSI}) holds with $M[\psi_n]\equiv 0$.

Otherwise let $c_\nu$ be the first non-zero coefficient of $\psi_n$ with index $\nu\ge 1$. Then
$$\begin{array}{rcl}\psi_n^{(k)}(z)&=&O(|z|^{-\nu-k}),\cr
w^{(k)}(z)&=&\psi^{(k)}_n(z)+o(|z|^{-n-k})=O(|z|^{-\nu-k}),\cr
(w^{(k)}(z))^{\ell_k}&=&(\psi^{(k)}_n(z))^{\ell_k}+O(|z|^{-(\nu+k)(\ell_k-1)})o(|z|^{-n-k})\cr
&=&(\psi^{(k)}_n(z))^{\ell_k}+o(|z|^{-k\ell_k-n-\nu(\ell_k-1)})\end{array}$$
and $M[w]=M[\psi_n]+R_n$ holds with remainder term
$$\begin{array}{rcl}R_n&=&\ds\sum_{\ell_j>0}\prod_{k\ne j}(\psi_n^{(k)}(z))^{\ell_k}\,o(|z|^{-j\ell_j-n-\nu(\ell_j-1)})\cr
&=&\ds\sum_{\ell_j>0}o\Big(|z|^{-\sum\limits_{k\ne j}\ell_k(k+\nu)-j\ell_j-n-\nu(\ell_j-1)}\Big)=o(|z|^{-d_M-n+1}),
\end{array}$$
since $\nu-\sum_{k=1}^m\ell_k(\nu-1)\le 1$, hence
$$\begin{array}{rcl}w^{\ell_0}M[w]&=&\ds\psi_n^{\ell_0}M[\psi_n]+o(|z|^{-n})M[\psi_n]
+o(|z|^{-d_M-n+1})\cr
&=&\ds\psi_n^{\ell_0}M[\psi_n]+
o\Big(|z|^{-n-\sum\limits_{k=1}^m(k+\nu)\ell_k}\Big)+o(|z|^{-d_M-n+1})\cr
&=&\ds\psi_n^{\ell_0}M[\psi_n]+o(|z|^{-d_M-n+1}).\quad\ende\end{array}$$

\subsection{\it Existence of solutions with specific asymptotic expansions} Theorem~\ref{AsympEx} says that some {\it specific} solution has an asymptotic expansion.
This may not be mixed up with the well-known and in some sense much more general Theorems~12.1.
and 14.1.\ in {\sc Wasow}~\cite{Wasow1}, which assert the existence of {\it some} solution having an asymptotic expansion. To the convenience of the reader we will prove existence of first, second, and fourth {\sc Painlev\'e} transcendents with prescribed asymptotic expansions
according to (\ref{asympexw}) on arbitrary sectors $S$ with central angle $\Theta=\frac45\pi$, $\Theta=\frac23\pi$, and $\Theta=\frac12\pi$, respectively.

\begin{thm}\label{EXISTENCE}To every such sector $S$  there exists some first, second, and fourth
{\sc Painlev\'e} transcendent with prescribed asymptotic expansion {\rm (\ref{asympexw})(i)}, {\rm (ii)}, and  {\rm (iv)} on $S$,
respectively.\end{thm}

\beweis According to our philosophy we prefer to consider equation (IV) in detail. We set $t=z^2$ and  $w(z)=t^{-1/2} v(t)^2$  to obtain
$$\ddot v=\frac{v^4-\gamma^2}{4v^3}+\frac1{4t}(2v^3-\alpha v)+\frac3{16t^2}(v^5-v)=f(t,v)\quad(\gamma^2=-\beta/2\ne 0),$$
hence
$\dot\x=\f(t,\x)$ when written as a system, with $\x=\left(\!\!\begin{array}{c} x\cr y\cr\end{array}\!\!\right),$
$x=v-\sqrt{\gamma},$ $y=\dot v$, $\f(t,\x)=\left(\!\!\begin{array}{c}y\cr f(t,x+\sqrt{\pm\gamma})\cr\end{array}\!\!\right),$
and Jacobian $\ds\lim_{t\to\infty}\f'(t,\x)=\left(\!\!\begin{array}{cc}0&1\cr 1&0\cr\end{array}\!\!\right)$.
From the latter it follows that there exists a unique formal solution $\sum\limits_{k=1}^\infty\x_kt^{-k}$, hence
a solution  with asymptotic expansion $\x(t)\sim \sum\limits_{k=1}^\infty\x_kt^{-k}$ on any given sector of central angle $\pi$, which itself
gives rise to a solution to (IV) with asymptotic expansion
(\ref{asympexw})(iv)$_c^\pm$ on any given sector of central angle $\pi/2$.
Similarly, the substitution $t=z^2$, $w(z)=t^{1/2}v(t)^2$ yields
$$\ddot v+\frac{\dot v}t=\frac1{16}v(v^2+2)(3v^2+2)-\frac{\alpha v}{4t}+\frac1{16t^2}\Big(v-\frac{4\gamma^2}{v^3}\Big),$$
hence we obtain solutions having asymptotic expansions {\rm (iv)}$_a$ and {\rm (iv)}$_b$ (according to $\lim_{t\to\infty}v(t)^2=-2/3$ and $\lim_{t\to\infty}v(t)^2=-2$), respectively. Again the sector $S$ in the $z$-plane with central angle $\pi/2$ may be prescribed.
In the other cases we just note the substitutions and differential equations:

\begin{tabular}{lrl}
(i)&  $z=t^{4/5}$, $w(z)=t^{2/5}v(t)$:& $\ds\ddot v+\frac{\dot v}t=\frac{96}{25}(6v^2+1)+\frac{4v}{25t^2}.$\cr
(ii)$_a$& $z=t^{2/3}$, $w(z)=t^{-2/3}v(t)$:& $\ds\ddot v-\frac{\dot v}t=\frac 49(v+\alpha)+\frac{8}{9t^2}(v^3-v).$\cr
(ii)$_b$& $z=t^{2/3}$,  $w(z)=t^{1/3}v(t)$:& $\ds\ddot v+\frac{\dot v}t=\frac 49(v+2v^3)+\frac{4\alpha}{9t}+\frac{v}{9t^2}.$ \Ende
\end{tabular}

\end{document}